\documentclass[12pt,a4]{amsart}
\usepackage{amsmath}
\usepackage{amsthm}
\usepackage{amssymb}
\usepackage{graphicx}
\usepackage{multirow}
\usepackage{lscape}
\usepackage[all]{xy}
\usepackage[thinlines,thiklines]{easybmat}

\newtheorem{rem}{Bemerkung}

\newcommand{\G}{{\mathbb G}}
\newcommand{\Z}{{\mathbb Z}}
\newcommand{\F}{{\mathbb F}}
\newcommand{\N}{{\mathbb N}}

\newcommand{\A}{{\mathbb A}}
\newcommand{\BQ}{{\mathbb Q}}
 \newcommand{\MF}{{\mathcal F}}
\newcommand{\MK}{{\mathcal K}}
\newcommand{\MT}{{\mathcal T}}
\newcommand{\MO}{{\mathcal O}}

\newcommand{\MC}{{\mathcal C}}

\newcommand{\BP}{{\mathbb P}}
\newcommand{\fp}{{\frak p}}
\newcommand{\La}{{\Lambda}}

\newcommand{\la}{{\lambda}}
\newcommand{\Ga}{{\Gamma}}
\newcommand{\ga}{{\gamma}}
\newcommand{\bo}{{\boldsymbol{\omega}}}

\newcommand{\fa}{{\boldsymbol{a}}}

\newcommand{\fg}{{\boldsymbol{g}}}
\newcommand{\fu}{{\boldsymbol{u}}}
\newcommand{\bn}{{\boldsymbol{n}}}
\newcommand{\fv}{{\boldsymbol{v}}}

\newcommand{\Om}{{\Omega}}

 \newcommand{\MB}{{\mathcal B}}
\newcommand{\lra}{{\longrightarrow}}

\newcommand{\om}{{\omega}}

\newcommand{\BF}{{\boldsymbol{F}}}

\makeindex
\numberwithin{equation}{section}

\makeatletter
\newcommand{\xyR}[1]{%
\makeatletter
\xydef@\xymatrixrowsep@{#1}
\makeatother
} 
\newcommand{\xyC}[1]{%
\makeatletter
\xydef@\xymatrixcolsep@{#1}
\makeatother
} 

\begin{document}

\setcounter{section}{-1}
\parindent=0pt

\title[On Drinfld modular forms of higher rank IV]
{On Drinfeld modular forms of higher rank IV: Modular forms with level} 
 \author{Ernst-Ulrich Gekeler}
 \maketitle
 \begin{abstract}
 We construct and study a natural compactification $\overline{M}^r(N)$ of the moduli scheme $M^r(N)$ for
 rank-$r$ Drinfeld $\F_q[T]$-modules with a structure of level $N \in \F_q[T]$. Namely, $\overline{M}^r(N) =
 {\rm Proj}\,{\bf Eis}(N)$, the projective variety associated with the graded ring ${\bf Eis}(N)$ generated by the Eisenstein
 series of rank $r$ and level $N$. We use this to define the ring ${\bf Mod}(N)$ of all modular forms of rank $r$ and level $N$. It
 equals the integral closure of ${\bf Eis}(N)$ in their common quotient field $\widetilde{\MF}_r(N)$. Modular forms are
 characterized as those holomorphic functions on the Drinfeld space $\Om^r$ with the right transformation behavior under
 the congruence subgroup $\Ga(N)$ of $\Ga = {\rm GL}(r,\F_q[T])$ (``weak modular forms'') which, along with all their
 conjugates under $\Ga/\Ga(N)$, are bounded on the natural fundamental domain $\BF$ for $\Ga$ on $\Om^r$.
 \end{abstract}
{\bf 0. Introduction.}
 \vspace{0.3cm}\\ 
(0.1) This is the fourth of a series of papers (see \cite{19}, \cite{20}, \cite{21}) which aim to lay the foundations for
a theory of Drinfeld modular forms of higher rank. These are modular forms for the modular group $\Ga = {\rm GL}(r,\F_q[T])$
or its congruence subgroups, where ``higher rank'' refers to $r$ larger or equal to 2. The case of $r=2$, remarkably similar
in some aspects but rather different in others to the theory of classical elliptic modular forms for ${\rm SL}(2,\Z)$ or
its congruence subgroups, is meanwhile well-established and the subject of several hundred publications since about 1980.
 \bigskip
 
We leave aside to deal with more general Drinfeld coefficient rings $A$ than $A=\F_q[T]$, as the amount of technical
and notational efforts required would obscure the overall picture. The interested reader may consult \cite{15} to get
an impression of the complications that - even for $r=2$ - result from class numbers $h(A) > 1$ for general $A$.
 \bigskip
 
 (0.2) While we developed some of the theory of modular forms ``without level'' in \cite{19} and \cite{21} and focussed on
 the connection with the geometry of the Bruhat-Tits building in \cite{19} and \cite{20}, the current part IV is devoted
 to forms ``with level'', i.e., forms for congruence subgroups of $\Ga$.  Again we restrict to the most simple case of full
 congruence subgroups $\Ga(N)=\{\ga\in \Ga\,|\, \ga \equiv 1 (\bmod N)\}$ for $N \in A$. Finer arithmetic/geometric properties
 of modular forms (or varieties) for other congruence subgroups $\Ga' \supset \Ga(N)$ may be derived in the course of 
 the further development of the theory from those for $\Ga(N)$, by taking invariants (or quotients) of the finite
 group $\Ga'/\Ga(N)$.
  \bigskip
  
(0.3) Let us introduce a bit of notation: $\F = \F_q$ is the finite field with $q$ elements, $A = \F[T]$ the polynomial ring in
an indeterminate $T$, with quotient field $K = \F(T)$, and its completion $K_{\infty} = \F((T^{-1}))$ at infinity, and $C_{\infty}$
the completed algebraic closure of $K_{\infty}$. The Drinfeld symmetric space $\Om^r$ (where $r \geq 2$) is the complement
in $\mathbb P^{r-1}(C_{\infty})$ of the $K_{\infty}$-rational hyperplanes. The modular group $\Ga = {\rm GL}(r,A)$ acts
in the usual fashion on $\Om^r$, and we let $M^r(N)$ be the quotient analytic space $\Ga(N)\setminus \Om^r$ (which
is also the set of $C_{\infty}$-points of an affine variety labelled by the same symbol, and which is smooth if $N \in A$
is non-constant).
 \bigskip
 
The modular forms dealt with will be holomorphic functions on $\Om^r$ with certain additional properties; so the theory
is ``over $C_{\infty}$''; we will only briefly touch on questions of rationality. 
 \bigskip 
 
(0.4) Our approach is based on 
 \begin{itemize}
  \item[(i)] the use of the natural fundamental domain $\BF$ for $\Ga$ on $\Om^r$ introduced in \cite{18}; it relies on the
  notion of successive minimum basis (SMB) of an $A$-lattice in $C_{\infty}$. On $\BF$, one may perform explicit 
  calculations;
   \item[(ii)] a natural compactification $\overline{M}^r(N)$ of $M^r(N)$, the {\em Eisenstein compactification},
   whose construction is influenced by but different from Kapranov's in \cite{28}.
   \end{itemize}
 The obvious examples of modular forms-to-be for $\Ga(N)$ are the Eisenstein series of level $N$. They generate a
 graded $C_{\infty}$-algebra ${\bf Eis}(N)$ (generated in dimension 1 if $N$ is non-constant), and $\overline{M}^r(N)$
 will be the associated projective variety ${\rm Proj}({\bf Eis}(N))$, see Theorem 5.9. It is a closed subvariety of a certain
 projective space $\mathbb P^{c-1}$, where $c$ is the number of cusps of $\Ga(N)$ (Corollary 4.7, Theorem 5.9), and is
 therefore supplied with a natural very ample line bundle $\mathfrak{M}$. We define strong modular forms of weight
 $k$ for $\Ga(N)$ as sections of $\mathfrak M^{\otimes k}$, and thereby get the graded ring ${\bf Mod}^{\rm st}(N)$
 of strong modular forms, which encompasses ${\bf Eis}(N)$.
  \bigskip
  
(0.5) The Eisenstein compactification is natural and explicit, and has good functorial properties (see Remark 5.8; it is, e.g.,
compatible with level change), but unfortunately we presently cannot assure that it is normal. Correspondingly, strong
modular forms are integral over ${\bf Eis}(N)$ (and in fact over ${\bf Mod} = {\bf Mod}(1)$, the ring of modular forms of
type 0 for $\Ga(1) = \Ga$), but we don't know whether ${\bf Mod}^{\rm st}(N)$ is integrally closed. We {\em define} the
Satake compactification $M^r(N)^{\rm Sat}$ of $M^r(N)$ as the normalization of $\overline{M}^r(N)$ (as Kapranov does)
and a modular form of weight $k$ for $\Ga(N)$ as a section of the pull-back of $\mathfrak{M}^{\otimes k}$
to $M^r(N)^{\rm Sat}$. This yields the
graded ring ${\bf Mod}(N)$ of all modular forms. Hence we have inclusions 
 $${\bf Eis}(N) \subset {\bf Mod}^{\rm st}(N) \subset {\bf Mod}(N)\leqno{(0.6)}$$
of finitely generated graded integral $C_{\infty}$-algebras, where ${\bf Mod}(N)$ is the integral closure of ${\bf Eis}(N)$ in
their common quotient field $\widetilde{\mathcal{F}}_r(N)$. Elements of ${\bf Mod}(N)$ have a nice characterization given by
Theorem 7.9: A weak modular form $f$ of weight $k$ is modular if and only if $f$, together with all its conjugates
$f_{[\ga]_k}$ ($\ga \in \Ga/\Ga(N))$, is bounded on the fundamental domain $\BF$. Further ${\bf Eis}(N)$ has always finite
codimension in ${\bf Mod}^{\rm st}(N)$ (Corollary 7.11), while $\dim({\bf Mod}(N)/{\bf Mod}^{\rm st}(N))$ is either zero or infinite,
according to whether $M^r(N)^{\rm Sat}$ agrees with $\overline{M}^r(N)$ or not (Corollary 7.14). Except for some examples
presented in Section 8, where the two compactifications and also the three rings in (0.6) agree (these examples depend
crucially on work of Cornelissen \cite{8} and Pink-Schieder \cite{32}), we don't know what happens in general: more research
is needed! At least $M^r(N)^{\rm Sat}$ is not very far from $\overline{M}^r(N)$: the normalization map
 $$\nu:\: M^r(N)^{\rm Sat} \lra \overline{M}^r(N) $$
 is bijective on $C_{\infty}$-points
 (Corollary 7.6, see also Proposition 1.18 in \cite{28}), and is an isomorphism on the complement of a closed subvariety of
 codimension $\geq 2$ (Corollary 6.10).
  \bigskip
  
(0.7) We now describe the plan of the paper. In the first section, we introduce the space $\overline{\Om}^r$ with its strong
topology, which upon dividing out the action of $\Ga(N)$ will yield the underlying topological space for the Eisenstein
compactification $\overline{M}^r(N)$. Its points correspond  to homothety classes of pairs $(U,i)$, where $U \not= 0$ is
a $K$-subspace of $K^r$ and $i$ is a discrete embedding of $U\cap A^r$ into $C_{\infty}$. For technical purposes we also 
consider the $\G_m$-torsor $\overline{\Psi}^r$ over $\overline{\Om}^r$ whose points correspond to pairs (not
homothety classes) $(U,i)$ as above. Further, the fundamental domains $\widetilde{\BF}$ on $\Psi^r$ and $\BF$ on
$\Om^r$ for $\Ga$ are introduced. Although $\overline{\Om}^r$ and $\overline{\Psi}^r$ come with the same information,
it will sometimes be more convenient to work with $\overline{\Psi}^r$ and $\widetilde{\BF}$ instead of $\overline{\Om}^r$
and $\BF$. We take particular care to give a consistent description of the group actions on $\overline{\Om}^r$ and related
objects.
 \bigskip
 
In Section 2 the (well-known) relationship of $\Om^r$ with the moduli of Drinfeld modules of rank $r$ is presented. We
further show the crucial technical result Theorem 2.3, which asserts that the bijection 
 $$j:\, \Ga \setminus \overline{\Om}^r\stackrel{\cong}{\lra} {\rm Proj}({\bf Mod})$$
is a homeomorphism for the strong topologies on both sides. We further 
introduce and describe the function fields of the analytic spaces $M^r(N) = \Ga(N)\setminus \Om^r$ and 
$\widetilde{M}^r(N) = \Ga(N)\setminus \Psi^r$.
 \bigskip
 
In Sections 3 and 4, the boundary components and the (non-)vanishing of Eisenstein series on them are studied. 
We find in Corollary 4.7 that the space ${\rm Eis}_k(N)$ of Eisenstein series of level $N$ and weight $k$ has dimension
 $c_r(N)$, the number of cuspidal divisors of $\Ga(N)\setminus \overline{\Om}^r$, independently of $k$. Further
 (Proposition 4.8), ${\rm Eis}_1(N)$ separates points of $\Ga(N)\setminus\overline{\Om}^r$, which will give rise to
 its projective embedding. This latter is defined and investigated in Section 5; we thereby interpret $\Ga(N)\setminus
 \overline{\Om}^r$ as the Eisenstein compactification $\overline{M}^r(N)$ of $M^r(N)$.
  \bigskip
  
 Section 6 is of a more technical nature. There we construct  tubular neighborhoods along the cuspidal divisiors of 
 $\overline{M}^r(N)$, see Theorem 6.9.
  \bigskip
  
In Section 7, the rings ${\bf Mod}^{\rm st}(N)$ and ¢${\bf Mod}(N)$ of modular forms are introduced and their relation with
the Eisenstein ring ${\bf Eis}(N)$ and the compactifications $\overline{M}^r(N)$ and $M^r(N)^{\rm Sat}$ is  discussed.
 \bigskip
 
We conclude in Section 8 with the two classes of examples where our knowledge is more satisfactory than in the
general situation, namely the special cases where either the rank $r$ equals $2$ or where the conductor $N$ has
degree 1.
 \bigskip
 
(0.8) The point of view (and the notation, see below) of this paper widely agrees with that of the preceding
\cite{19}, \cite{20}, \cite{21}, to which we often refer. As in these, our basic references for rigid analytic geometry 
are the books \cite{12} by Fresnel-van der Put and \cite{7} of Bosch-G\"untzer-Remmert. The canonical topology on the
set $X(C_{\infty})$ of $C_{\infty}$-points of an analytic space $X$ (\cite{7} Section 7.2) is labelled as the
{\em strong topology}, so functions continuous with respect to it are strongly continuous, etc. In general, we don't
distinguish in notation between $X$ and $X(C_{\infty})$; ditto, a $C_{\infty}$-variety and its analytification are usually
described by the same symbol. It will (hopefully) always be clear from the context whether e.g. the ``algebraic'' or the
``analytic'' local ring is intended. 
 \bigskip
 
(0.9) After this paper was largely completed, I got access to the recent preprints \cite{4}, \cite{5}, \cite{6} of Dirk
Basson, Florian Breuer, and Richard Pink, which go about the same topic: providing a foundation for the theory
of higher rank Drinfeld modular forms. As it turns out, the relative perspectives of Basson-Breuer-Pink's work and of
the current paper are rather different. While BBP deal with the most general Drinfeld coefficient rings $A$ and arithmetic
subgroups of ${\rm GL}(r,A)$, for which they establish basic but sophisticated facts like e.g. the existence of expansions 
around infinity of weak modular forms, we restricted to the coefficient ring $A=\F_q[T]$ and full congruence subgroups
and focus on the role of Eisenstein series, their arithmetic properties, and their impact on compactifications of the moduli
schemes. Apart from examples, there is little overlap between the two works; so the reader who wants to enter into the 
field might profit from studying the two of them.
 \medskip
 
Finally, I wish to point to the recent thesis \cite{245} of Simon H\"aberli, whose purpose is similar. In contrast with the
present article, H\"aberli gives a direct construction of the Satake compactification, which he uses for the description of 
modular forms. 
 \medskip

{\bf Notation.}
\vspace{0.3cm}

$\F = \F_q$ the finite field with $q$ elements;\\
$A=\F[T]$ the polynomial ring in an indeterminate $T$, with quotient field $K=\F(T)$ and its completion $K_{\infty}=
\F((T^{-1}))$ at infinity;\\
$C_{\infty}=$ completed algebraic closure of $K_{\infty}$, with absolute value $|~.~|$ and valuation 
$v:\: C_{\infty}^{\ast} \lra \BQ$ normalized by $v(T) = -1$, $|T| = q$;\\
$\Psi^r = \{\bo = (\om_1,\ldots,\om_r) \in C_{\infty}^r~|~\mbox{the $\om_i$ are $K_{\infty}$-linearly independent}\}$\\
$\Om^r = \{\bo = (\om_1:\ldots: \om_r) \in \mathbb P^{r-1} (C_{\infty})~|~\bo \mbox{ represented by } (\om_1,\ldots,\om_r)
\in \Psi^r\}$\\
$\Ga = \Ga_r ={\rm GL}(r,A)$ with center $Z \cong \F^{\ast}$ of scalar matrices;\\
$\Ga(N) = \{\ga \in \Ga~|~\ga \equiv 1 (\bmod N)\}$, $N \in A$;\\
$\mathfrak{U} = \mbox{set of $K$-subspaces $U \not= 0$ of $V = K^n$}$;\\
$\Psi_U \cong \Psi^s$, $\Om_U \cong \Om^s$ attached to $U \in \mathfrak U$, where $\dim U = s$;\\
$\overline{\Psi}^r  = \underset{U\in \mathfrak U}{\stackrel{\bullet}{\bigcup}}\Psi_U$, 
$\overline{\Om}^r  = \underset{U\in \mathfrak U}{\stackrel{\bullet}{\bigcup}} \Om_U.$
\smallskip
 
If the group $G$ acts on the space $X$ then $G_x$, $Gx$ and $G\setminus X$ denote the stabilizer of $x \in X$, its orbit,
and the space of all orbits, respectively. Also, for $Y \subset X$, $G\setminus Y$ is the image of $Y$ in $G\setminus X$. The 
multiplicative group of the ring $R$ is $R^{\ast}$; the $R$-module generated by $x_1,\ldots,x_r$ is written either 
as $\sum Rx_i$ or as $\langle x_1,\ldots,x_r\rangle_R$. We use the convention $\N=\{1,2,3,\ldots\} $, $\N_0 =\{0,1,2,\ldots\}$.
 \vspace{0.5cm}
 
 {\bf 1. The spaces $\overline{\Psi}^r$ and $\overline{\Om}^r$.}
  \vspace{0.3cm}
  
  (1.1) We let $V$ be the $K$-vector space $K^r$, where $r \geq 2$, and $\mathfrak U$ the set of $K$-subspaces $U \not= 0$ of $V$. An $A$-lattice in $U \in \mathfrak U$ is a free $A$-submodule $L$ of $U$ of full rank ${\rm rk}_A(L) = \dim_K(U)$, that is $K \otimes
  L = KL = U$. A subset of $C_{\infty}$ is {\em discrete} if the intersection with each ball of finite radius in $C_{\infty}$ is
  finite. A {\em discrete embedding} of $U \in \mathfrak U$ (``embedding'' for short) is some $K$-linear injective map
  $i:\: U \lra C_{\infty}$ such that $i(L)$ is discrete in $C_{\infty}$ for one fixed (or equivalently, for each) $A$-lattice $L$ in $U$.
  We put
   $$\begin{array}{lll}
   \Psi_U & := & \mbox{set of discrete embeddings of $U$, and}\smallskip\\
    \Om_U & := & C_{\infty}^{\ast}\setminus \Psi_U, \mbox{the quotient of $\Psi_U$ modulo the}\\
     & & \mbox{action of the multiplicative group $C_{\infty}^{\ast}$.}
      \end{array}\leqno{(1.2)}$$
Further, $\Psi^r := \Psi_V = \Psi_{K^r}$, $\Om^r := \Om_V$, and 
 $$\overline{\Psi}^r := \underset{U\in \mathfrak U}{\stackrel{\bullet}{\bigcup}} \Psi_U,\: \overline{\Om}^r :=
  \underset{U\in \mathfrak U}{\stackrel{\bullet}{\bigcup}} \Om_U.$$  
If $U \subset U' \in \mathfrak A$, restriction to $U$ defines canonical maps
 $$\Psi_{U'} \lra \Psi_U \mbox{ and } \Om_{U'} \lra \Om_U. \leqno{(1.2.1)}$$
(1.3) We let $L_V := A^r$ and $L_U := L_V \cap U$ be the standard lattices in $V$ and $U$, respectively. As a $K$-linear map
$i:\: V \lra C_{\infty}$ is discrete if and only if the images $\om_j := i(e_j)$ of the standard basis vectors $e_j$
($1 \leq j \leq r$) are $K_{\infty}$-linearly independent (l.i.), we see that
 $$\Psi^r = \{\bo = (\om_1,\ldots,\om_r) \in C_{\infty}^{\ast}~|~\om_1,\ldots,\om_r \mbox{ l.i.}\}.$$
After choosing bases of the subspaces $U$, we get similar descriptions for $\Psi_U$ and the quotients $\Om^r$ and $\Om_U$. In
particular, we find for $r=2$ the familiar Drinfeld upper half-plane
 $$\begin{array}{rll}
 \Om^2 = C_{\infty}^{\ast}\setminus \{(\om_1,\omega_2)~|~\om_1,\om_2 \mbox{ l.i.}\} &\stackrel{\cong}{\lra}& C_{\infty} 
 \setminus K_{\infty}.\\
  (\om_1,\om_2) & \longmapsto & \om_1/\omega_2 \end{array}$$
  (1.4) The sets $\Psi^r$ and $\Om^r$ (and therefore also $\Psi_U$ and $\Om_U$) are equipped with structures of
 $C_{\infty}$-analytic spaces (actually defined over $K_{\infty}$), namely as admissible open subspaces of $\A^r(C_{\infty}) =
 C_{\infty}^r$ or of $\mathbb P^{r-1}(C_{\infty})$, respectively, see \cite{11},\cite{10}, or \cite{33}.
  \bigskip
  
(1.5) The group ${\rm GL}(r,K)$ acts as a matrix group from the right on $V$, which induces left actions on 
$\overline{\Psi}^r$ and $\overline{\Om}^r$, viz.: For $\ga \in {\rm GL}(r,K)$, let $r_{\ga}:\: V \lra V$ be the map 
$x \longmapsto x\ga$. Then $\ga$ maps $(U,i:\: U \hookrightarrow C_{\infty}) \in \Psi_U$ to $\ga(U,i) := (U\ga^{-1},
i \circ r_{\ga})$. The reader may verify that this, together with the description of $\Psi^r$ in (1.3), yields the standard
left matrix action of $\ga$ on $\Psi^r$, the elements of $\Psi^r$ being regarded as column vectors 
$(\om_1,\ldots,\om_r)^t$.
 \bigskip
 
(1.6) Since $A$ is a principal ideal domain, the theory of finitely generated modules over such (e.g. \cite{29} XV Sect. 2)
shows that $\Ga := {\rm GL}(r,A)$ acts transitively on the set $\mathfrak U_s$ of $U \in \mathfrak U$ of fixed dimension 
$s$. We use as a standard representative for $\mathfrak U_s$ the space
 $$V_s := \{(0,\ldots,0,*,\ldots,*) \in V\}\leqno{(1.6.1)}$$
 of vectors whose first $r-s$ entries vanish. The fixed group of $V_s$ ($1 \leq s <r$) in ${\rm GL}(r,K)$ is the maximal
 parabolic subgroup 
  $$P_s := \left\{
   \begin{array}{|c|c|}\hline
   * & *\\ \hline
   0 & *\\ \hline \end{array} \right\}
   \leqno{(1.6.2)}$$
of matrices with an $(r-s,s)$-block structure whose lower left block vanishes. The action of $P_s$ on $V_s$ is via the
group
$$M_s := \left\{
   \begin{array}{|c|c|}\hline
   1 & 0\\ \hline
   0 & *\\ \hline \end{array} \right\}
   \leqno{(1.6.3)}$$
regarded as a factor group of $P_s$.	
 \bigskip
 
(1.7) As explained in (1.3), the choice of a $K$-basis of $U \in \mathfrak U$ yields an embedding of $\Psi_U$ into 
$C_{\infty}^{\dim U}$. The Haussdorff topology induced on $\Psi_U$ is independent of that choice, and is referred to as
the {\em strong topology} on $\Psi_U$. Similarly, using embeddings into projective spaces, we define the strong 
topologies on the $\Om_U$.
 \bigskip
 
(1.8) Our next aim is to define reasonable strong topologies on $\overline{\Psi}^r$ and $\overline{\Om}^r$ extending 
the topologies on the strata. For this we recall the concept of successive minimum bases. An {\em $A$-lattice in}
$C_{\infty}$ is a discrete $A$-submodule $\La$ of finite rank. A {\em successive minimum basis} (SMB) of 
$\La$ is an ordered $A$-basis $\{\om_1,\ldots,\om_r\}$ of $\La$ (note this differs from usual set-theoretic notation)
subject to: For each $1 \leq j \leq r$, $|\om_j|$ is minimal among
 $$\{|\om|~|~\om \in \La \setminus (A\om_1+ \cdots + A\om_{j-1})\}.$$
(For $j = 1$ this means: $\om_1$ is a lattice vector of minimal non-zero length.) It is shown in \cite{18} Proposition 3.1 that
each $A$-lattice $\La$ in $C_{\infty}$ possesses an SMB $\{\om_1,\ldots,\om_r\}$, and it has the following  additional
properties:
 \bigskip
 
(1.8.1) The $\om_i$ are orthogonal, that is, given $a_1,\ldots,a_r \in K_{\infty}$,
 $$|\sum_{1 \leq i \leq r} a_i\om_i| = \underset{i}{\max} |a_i||\om_i|;$$    
(1.8.2) The series of positive real numbers $|\om_1| \leq |\om_2| \leq \ldots \leq |\om_r|$ is an invariant of $\La$, that is,
independent of the choice of the SMB.
 \bigskip
 
(1.9) We define the strong topology on $\overline{\Psi}^r$ as the unique Hausdorff topology which satisfies for each
$U' \in \mathfrak A$:
  \bigskip
  
 (1.9.1) Restricted to $\Psi_{U'}$, it agrees with the strong  topology given there by (1.7);
  \bigskip
  
(1.9.2) The topological closure $\overline{\Psi}_{U'}$ of $\Psi_{U'}$ equals $\underset{U \subset U'}{\stackrel{\bullet}{\bigcup}}\Psi_U$;
 \bigskip
 
(1.9.3) Assume $U' \supset U \in \mathfrak U$, and let $i:\: U \lra C_{\infty}$ and $i_k:\: U' \lra C_{\infty}$ ($k \in \N$) be
discrete embeddings. Then $(U,i) = {\displaystyle \lim_{k\to \infty} (U',i_k)}$ if and only if     
 \begin{itemize}
  \item[(a)]for each $\la \in L_U$, $i(\la) = {\displaystyle \lim_{k\to \infty} i_k(\la)}$ and 
   \item[(b)] for each $\la \in L_{U'} \setminus L_U$, ${\displaystyle \lim_{k\to \infty} |i_k(\la)| = \infty}$, uniformly in $\la$.
  \end{itemize}
 Note that it suffices to require (a) for the elements of a basis of $L_U$. In qualitative terms, $i_k:\: U' \hookrightarrow
 C_{\infty}$ is very close to $i:\: U \hookrightarrow C_{\infty}$ iff
  \begin{itemize}
  \item[(a')] $i_k(\la)$ is very close to $i(\la)$ for the elements $\la$ of an $A$-basis of $L_U$, and
  \item[(b')] for each $\la \in L_{U'} \setminus L_U$, $|i_k(\la)|$ is very large compared to the $|\om_j|$, where $\{\om_j\}$
  is an SMB of $i(L_U)$.
  \end{itemize}
  Furthermore, (b') my be replaced by
   \begin{itemize}
   \item[(b'')] for each $\la \in L_{U'}\setminus L_U$, $i_k(\la)$ has very large distance $d(i_k(\la),K_{\infty}i(U))$ to the
 $K_{\infty}$-space generated by $i(U)$.
\end{itemize}  
The strong topology on $\overline{\Om}^r = C_{\infty}^{\ast} \setminus  \overline{\Psi}^r$ is the quotient topology; it has
properties analogous to (1.9.1)--(1.9.3). Obviously, the action of ${\rm GL}(r,K)$ on both $\overline{\Psi}^r$ and 
$\overline{\Om}^r$ is through homeomorphisms w.r.t. the so defined topologies.
 \bigskip
 
(1.10) A continuous function $f:\: \overline{\Psi}^r \lra C_{\infty}$ has {\em weight} $k \in \Z$ if 
 $$f(U,c \cdot i) = c^{-k} f(U,i)$$
 holds for $c \in C_{\infty}^{\ast}$ and $(U,i) \in \overline{\Psi}^r$.
  \bigskip

(1.11)  The basic examples of functions with weight are the various types of {\em Eisenstein series} defined below.
For $k \in \N$ put
 $$E_k(U,i) := \underset{\la \in L_U}{\sum{'}} \,i(\la)^{-k}.$$
(The prime $\sum'$ indicates that the sum is over the non-zero elements of the index set.) The following are
obvious or easy to show:
 \begin{itemize}
 \item[(i)] The sum converges and defines a continuous (even analytic) function $E_k$ on $\Psi_U$ ($U \in \mathfrak A$),
 which is non-trivial if and only if $k \equiv 0 (\bmod\, q-1)$;
  \item[(ii)] $E_k$ is continuous on the whole of $\overline{\Psi}^r$ with respect to the strong topology (due to the very
  definition of the latter);
   \item[(ii)] $E_k$ has weight $k$;
    \item[(iv)] $E_k$ is invariant under ${\rm GL}(r,A)$. 
 \end{itemize}
  (1.12) Now let $N$ be a non-constant monic element of $A$ and $\Ga(N) = \{\ga \in \Ga~|~\ga \equiv 1 (\bmod N)\}$ be the
  {\em full congruence subgroup} of level $N$. Fix some vector $\fu = (u_1,\ldots,u_r) \in V = K^r$ with $N\fu \in L_V = A^r$,
  and put
   $$E_{k,\fu}(U,i) := \underset{\la \in U \atop \la \equiv \fu(\bmod L_V)}{\sum^{}{'}}i(\la)^{-k}.$$
The following hold:\\
(i') The sum converges and defines a continuous (even analytic) function $E_{k,\fu}$ on $\Psi_U$; it depends only on the
residue class of $\fu$ modulo $L_V$, and is called the {\em partial Eisenstein series} with congruence condition $\fu$;\\
(ii), (iii) (see (1.11)), and\\
(iv') $E_{k,\fu\ga}(U,i) = E_{k,\fu}(\ga(U,i))$, $\ga \in \Ga$. In particular, $E_{k,\fu}$ is invariant under $\Ga(N)$.
 \bigskip
 
{\bf Remark.} 
 $$E_{k,\fu}(U,i) = N^{-k}\underset{\la \in L_U \atop \la \equiv N\fu(\bmod NL_V)}{\sum{'}} i(\la)^{-k},$$
 which up to the factor $N^{-k}$ is a partial sum of $E_k(U,i)$. This explains the notation ``partial Eisenstein series''.
  \bigskip
  
(1.13) By definition, $\om_r \not= 0$ for $\bo = (\om_1,\ldots,\om_r) \in \Psi^r$. Therefore we can normalize projective
coordinates on $\Om^r \subset \mathbb P^{r-1}(C_{\infty})$ so that
 \medskip
 
(1.13.1) $ \om_r=1$, i.e.,
 $$
 \Om^r =  \{(\om_1,\ldots,\om_{r-1})=(\om_1:\ldots:\om_{r-1}:1)  ~|~
  \om_1,\ldots,\om_{r-1},\om_r=1 \mbox{ l.i.}\}. 
 $$
 Similarly we usually assume $\om_r=1$ for $\bo = (\om_1:\ldots:\om_r) \in \Om_U$ if $U$ is one of the spaces
 $V_s$ of (1.6). With that convention, the Eisenstein series $E_k$ and $E_{k,\fu}$ may be regarded as functions on
 $\underset{1 \leq s \leq r}{\bigcup} \Om_{V_s}$. If $\bo \in \Om^r$ then (iii), (iv), (iv') imply
  \bigskip
 $$E_k(\ga\bo) = {\rm aut}(\ga,\bo)^k E_k(\bo) \leqno{(1.13.2)}$$   
 and
 $$E_{k,\fu}(\ga\bo) = {\rm aut}(\ga,\bo)^k E_{k,\fu\ga}(\bo).\leqno{(1.13.3)}$$
Here $\ga \in \Ga$, $\bo = (\om_1:\ldots:\om_r)$ with $\om_r=1$, and ${\rm aut}(\ga,\bo)$ is the factor of automorphy
 $${\rm aut}(\ga,\bo) = \sum_{1 \leq i \leq r}\ga_{r,i} \om_i \not= 0.\leqno{(1.13.4)}$$
 We assign no value to $E_k(\bo)$ or $E_{k,\fu}(\bo)$ if $\bo = C_{\infty}^*(U,i) \in \overline{\Om}^r$ does not belong
 to $\bigcup \Om_{V_s}$, but are content with the distinction (always well-defined) of whether $E_k$ (resp.
 $E_{k,\fu}$) vanishes at $\bo$ or not.
  \bigskip
  
{\bf 1.14 Remark} (on notation). In order to avoid notational overflow, we use the same symbol $E_k$ for both
occurrences: as a $\Ga$-invariant function on $\overline{\Psi}^r$ of weight $k$, or as a function on $\bigcup  
\Om_{V_s}$ subject to (1.13.2). A similar remark applies to $E_{k,\fu}$ and to other functions with weight.
 \bigskip

(1.15) We finally define fundamental domains for the actions of $\Ga$ on $\Psi^r$ and $\Om^r$. To wit, put
 $$\begin{array}{lll}
  \widetilde{\BF} & := &\{\bo = (\om_1,\ldots,\om_r)\in \Psi^r~|~\{\om_r,\om_{r-1},\ldots,\om_1\} \mbox{ is an SMB}\\
  &&\hspace*{4cm} \mbox{of its lattice } \La_{\bo} = \langle \om_1,\ldots,\om_r\rangle_A\}\\
   \BF &:=& C_{\infty}^*\setminus \tilde{\BF}. \end{array}$$
 (Note the reverse order of the $\om_i$!) They have the following properties.
    \bigskip
    
(1.15.1) As the condition for $\bo \in \widetilde{\BF}$ is stable under the multiplicative group, $\widetilde{\BF}$ is the full cone
above $\BF$. 
 \bigskip
 
(1.15.2) Each $\bo \in \Psi^r$ (resp. $\Om^r$) is $\Ga$-equivalent with at least one and at most a finite number of
$\bo' \in \widetilde{\BF}$ (resp. $\bo' \in \BF$).
 \bigskip
 
 \begin{proof}
 It suffices to treat the case $\widetilde{\BF}$. As each $A$-lattice $\La$ in $C_{\infty}$ has an SMB, the existence of a 
 representative $\bo' \in \tilde{\BF}$ for $\bo \in \Psi^r$ is obvious. Given $\bo \in \widetilde{\BF}$, the condition 
 $\ga\bo \in \widetilde{\BF}$ on $\ga \in \Ga$ together with (1.8.1) leads to bounds on the entries of $\ga$, which can
 be satisfied for a finite number of $\ga$'s only.
  \end{proof} 
  
 (1.15.3) $\widetilde{\BF}$ resp. $\BF$ is an admissible open subspace of $\Psi^r$ resp. $\Om^r$. 
  \bigskip
  
The most intuitive way to see this comes from identifying $\BF$ as the inverse image under the building map
$\la:\: \Om^r \lra \MB\MT(\BQ)$ of a subcomplex $W$ of the Bruhat-Tits building $\MB \MT$ of ${\rm PGL}(r,K_{\infty})$:
see \cite{19} Sect. 2. In fact, $W$ is a fundamental domain for $\Ga$ on $\MB\MT$.
 \bigskip
 
In view of the above, we refer to $\widetilde{\BF}$ resp. $\BF$ as the {\em fundamental domain} for $\Ga$ on $\Psi^r$
resp. $\Om^r$. As uniqueness of the representative in $\widetilde{\BF}$ resp. $\BF$ fails, this is weaker than the classical
notion of fundamental domain, but is still useful. Property (1.8.1) turns out particularly valuable for explicit calculations
with modular forms, as exemplified in \cite{19}. Also useful is the following observation, which is immediate from definitions. 
We formulate it for $\BF$ only, but it holds true also for $\widetilde{\BF}$.
 \bigskip
 
(1.15.4) Let $\BF_s$ be the fundamental domain for $\Ga_s = {\rm GL}(s,A)$ in $\Om_{V_s} \stackrel{\cong}{\lra} \Om^s$
($1 \leq s \leq r$). Then the strong closure of $\BF$ in $\overline{\Om}^r$ is $\overline{\BF} = \underset{1\leq s \leq r}{\bigcup}
F_s$. Each point of $\overline{\Om}^r$ is $\Ga$-equivalent with at least one and at most a finite number of points of
$\overline{\BF}$.
 \bigskip
 
Therefore, we can regard  $\overline{\BF}$ as a fundamental domain for $\Ga$ on $\overline{\Om}^r$. 
 \vspace{0.5cm}
 
{\bf 2. Quotients by congruence subgroups and moduli schemes.}
 \vspace{0.3cm} 
 
(2.1) Given an $A$-lattice $\La$ in $C_{\infty}$ of rank $r \in \N$, we dispose of \\
$\bullet$ the exponential function $e_\La:\: C_{\infty} \lra C_{\infty}$
 $$e_{\La}(z) = z\underset{\la \in \La}{\prod}'(1-z/\la) = \sum_{i\geq 0} \alpha_i(\La)z^{q^{i}}; \leqno{(2.1.1)}$$ 
$\bullet$ the Drinfeld $A$-module $\phi^{\La}$ of rank $r$, defined by  the operator polynomial
 $$ \phi_T^{\La}(X) = TX + g_1(\La) X^q+\cdots + g_r(\La)X^{q^r}, \mbox{ and}\leqno{(2.1.2)}$$
$\bullet$ the Eisenstein series
 $$E_k(\La) = \underset{\la \in \La}{\sum}' \la^{-k} \quad (k \in \N).\leqno{(2.1.3)}$$
We further put $g_0(\La) = T$, $E_0(\La) = -1$. These are connected by
 $$e_{\La}(Tz) = \phi_T^{\La}(e_{\La}(z));  \leqno{(2.1.4)}$$
 $$
 \underset{i,j \geq 0 \atop i+j =k}{\sum} \alpha_i E_{q^j-1}^{q^{i}} = \underset{i+j=k}{\sum} \alpha_i^{q^{j}} E_{q^{j}-1} = 1
   \mbox{ if } k = 0 \mbox{ and }
   0 \mbox{ otherwise},\leqno{(2.1.5)}$$
which determines a number of further relations, see e.g. \cite{16} Sect. 2.
 \bigskip
 
If $\La = \La_{\bo} = \underset{1 \leq i \leq r}{\sum} A\om_i$ with $\bo = (\om_1,\ldots,\om_r) \in \Psi^r$, we use  
$\bo$ instead of $\La$ as the argument. Thus $\phi^{\bo} = \phi^{\La_{\bo}}$, $e_{\bo} = e_{\La_{\bo}}$, etc. As functions
on $\Psi^r$, $g_i$, $\alpha_i$ are - like the Eisenstein series - holomorphic and $\Ga$-invariant of weight $q^{i}-1$, while
considered as functions on $\Om^r$, $g_i$ (and $\alpha_i$) satisfies
 $$g_i(\ga \bo) = {\rm aut}(\ga,\bo)^{q^{i}-1} g_i(\bo)$$
 (see Remark 1.14).
  \bigskip
  
(2.1.6) The three systems of functions on $\Psi^r$: $\{g_1,\ldots,g_r\}$, $\{\alpha_1,\ldots,\alpha_r\}$, $\{E_{q^{i}-1}~|~~
1 \leq i \leq r\}$ are each algebraically independent, and the relations between them are such that the ring
 $$ {\bf Mod} = \bigoplus_{k\leq 0} {\bf Mod}_k = C_{\infty}[g_1,\ldots,g_r],$$
 graded by the weight ${\rm wt}(g_i) := g^{i}-1$, may also be described as 
  $$C_{\infty}[\alpha_1,\ldots \alpha_r] = C_{\infty}[\alpha_i~|~i\in \N] = C_{\infty}[E_{q^{i}-1}~|~1 \leq i \leq r] =
   C_{\infty}[E_{q^{i}-1}~|~i\in \N].$$
(Actually {\bf Mod} is the ring of modular forms of type 0 for $\Ga$, see \cite{19}.)
 \bigskip
 
 (2.1.7) As a consequence, since the $g_i$ and $\alpha_i$ may be expressed through Eisenstein series, they have 
 strongly continuous extensions to $\overline{\Psi}^r$ and may therefore be evaluated on arbitrary points
 $\bo = (U,i) \in \overline{\Psi}^r$.
  \bigskip
  
(2.2) The Drinfeld modules $\phi^{\bo}$ and $\phi^{\bo'}$ ($\bo, \bo' \in \Om^r$) are isomorphic if and only if 
$\bo' = \ga\bo$ with some $\ga \in \Ga$. 
 \bigskip

Hence the map
 $$\begin{array}{rll}
  j:\: \Ga \setminus \Om^r & \hookrightarrow & {\rm Proj}\,{\bf Mod}\\
   \bo & \longmapsto & (g_1(\bo):\ldots: g_r(\bo))
  \end{array}$$
identifies the quotient analytic space of $\Om^r$ modulo $\Ga$ with the complement of the vanishing locus
of $\Delta := g_r$ in the weighted projective space $\overline{M}^r = {\rm Proj}\,{\bf Mod}$. (We remind the
reader that we do not distinguish in notation between a $C_{\infty}$-variety, its associated analytic space, and the
set of its $C_{\infty}$-points.) Here the $g_i$ are considered as formal variables of weight $q^{i}-1$, that is
$(x'_1:\ldots:x'_r) = (x_1:\ldots:x_r)$ in ${\rm Proj}\,{\bf Mod}$ if and only if there exists $c \in C_{\infty}^*$
such that $x'_i = c^{q^{i}-1} x_i$ for all $i$. In other words, via $j$
 $$\Ga \setminus \Om^r \stackrel{\cong}{\lra} M^r := ({\rm Proj}\,{\bf Mod})_{(g_r \not= 0)}\leqno{(2.2.1)}$$
equals (the set of $C_{\infty}$-points of) the moduli scheme $M^r$ for Drinfeld $A$-modules of rank $r$ over
$C_{\infty}$. The natural compactification of $M^r$ is
 $$ {\rm Proj}\,{\bf Mod} = \overline{M}^r = M^r \cup M^{r-1} \cup \cdots M^1, \leqno{(2.2.2)}$$
where for $1 \leq s \leq r$,
 $$(\Ga \cap P_s)\setminus \Om_{V_s} = {\rm GL}(s,A)\setminus \Om^s \stackrel{\cong}{\lra} M^s \quad 
  \mbox{(see (1.6))}$$
and $\Om^1 = M^1 = \{\mbox{point}\}$. Hence the stratification of the variety $\overline{M}^r$ corresponds to
that of    
 $$\Ga\setminus \overline{\Om}^r = \Ga \setminus (\underset{1 \leq s \leq r \atop U \in \mathfrak U_s}
  {\stackrel{\bullet}{\bigcup}} \Om_U) = \underset{1 \leq s \leq r}{\stackrel{\bullet}{\bigcup}} {\rm GL}(s,A)\setminus \Om^s
  \leqno{(2.2.3)}$$
under the bijection 
 $$\begin{array}{rll}
  j:\: \Ga \setminus \overline{\Om}^r & \stackrel{\cong}{\lra} & \overline{M}^r \\
   \bo & \longmapsto & (g_1(\bo): \ldots : g_r(\bo)),
  \end{array} \leqno{(2.2.4)}$$
which is well-defined in view of (2.1.7).
 \bigskip
 
In a similar way (although this looks a bit artificial), we may describe $\Ga\setminus \Psi^r$ via 
 
    $$\begin{array}{rcl}
  \tilde{j}:\: \Ga \setminus \Psi^r & \hookrightarrow & \A^r (C_{\infty})\\
   \bo& \longmapsto & (g_1(\bo), \ldots , g_r(\bo)),
  \end{array} \leqno{(2.2.5)}$$ 
as the complement $\widetilde{M}^r$ of ($g_r=0$) in $\A^r$. It is the moduli scheme of rank-$r$ Drinfeld $A$-modules 
over $C_{\infty}$  with a ``non-vanishing differential'', that is, with an identification of the underlying 
additive group with $\G_a$ or, what is the same, with explicit coefficients $g_i$ of its $T$-operator polynomial. The
horizontal compactification $\Ga \setminus \overline{\Psi}^r$ then becomes
 $$\begin{array}{lcl}
  \Ga \setminus \overline{\Psi}^r &= & \Ga \setminus (\underset{1 \leq s \leq r \atop U \in \mathfrak U_s}
  {\stackrel{\bullet}{\bigcup}} \Psi_U) = \underset{1 \leq s \leq r}{\stackrel{\bullet}{\bigcup}}
  {\rm GL}(s,A)\setminus \Psi^s\\
   & \underset{\widetilde{j}}{\stackrel{\cong}{\lra}} & \underset{1 \leq s \leq r}{\stackrel{\bullet}{\bigcup}}
    \widetilde{M}^s =: \overline{\widetilde{M}}^r = C_{\infty}^r \setminus\{0\}, \end{array} 
    \leqno{(2.2.6)}$$
in analogy with (2.2.2), (2.2.3), (2.2.4).
 \bigskip
 
(2.2.7) In the sequel, whenever writing $\Ga\setminus \Om^r = M^r$ or $\Ga\setminus \Psi^r = \widetilde{M}^r$, the
identification is via $j$ or $\widetilde{j}$, respectively.
 \bigskip
  
{\bf 2.3 Theorem.} {\it The map $j:\: \Ga \setminus \overline{\Om}^r \stackrel{\cong}{\lra} {\rm Proj}\,{\bf Mod} =
{\rm Proj}\,C_{\infty}[g_1,\ldots,g_r]$ of $(2.2.4)$ is a strong homeomorphism, i.e., with respect to the strong topologies 
on both sides. Similarly, $\widetilde{j}:\: \Ga \setminus\overline{\Psi}^r \stackrel{\cong}{\lra} \overline{\widetilde{M}}^r$
is a strong homeomorphism.}

 \begin{proof}    
The proof for $j$ will also show the statement for $\widetilde{j}$.
 \begin{itemize}
\item[(i)] By construction, $j$ is continuous as a map from $\overline{\Om}^r$, and thus as a map from $\Ga \setminus 
\overline{\Om}^r$ supplied with the quotient topology. Therefore we must show that $j^{-1}$ is continuous.
 \item[(ii)] Let $(\phi^{(n)})_{n\in \N}$ be a series of Drinfeld modules of rank $\leq r$, given by their $T$-division
 polynomials $\phi^{(n)}_T(X) = \underset{0 \leq i \leq r}{\sum} g_i^{(n)} X^{q^{i}}$ and converging to $\phi$ with
 $\phi_T(X) = \sum g_i X^{q^{i}}$. This means that $\fg^{(n)} = (g_1^{(n)}: \ldots : g_r^{(n)})$ converges to 
$\fg = (g_1:\ldots:g_r)$. Let $s$ be the rank of $\phi$, i.e., $g_s \not= 0$, $g_{s+1} = \cdots = g_r = 0$. We may
suppose that $g_s = \underset{n\to \infty}{\lim} g_s^{(n)} = 1$. Let $\La^{(n)}$ (resp. $\La$) be the lattice associated
to $\phi^{(n)}$ (resp. $\phi$), each provided with an SMB $\{\om_r^{(n)},\om_{r-1}^{(n)},\ldots,\om_1^{(n)}\}$
(resp. $\{\om_r,\ldots,\om_{r-s+1}\}$), where we have put $\om_i^{(n)} = 0$ for $i \leq r-{\rm rk}(\phi^{(n)})=r-{\rm rk}_A(\La^{(n)})$.
Put $\bo := (0:\ldots:0:\om_{r-s+1}:\ldots:\om_r)$ and $\bo^{(n)} := (\om_1^{(n)}: \cdots: \om_r^{(n)})$, and let
$[\bo]$ resp. $[\bo^{(n)}]$ be the corresponding class in $\Ga\setminus \overline{\Om}^r$. Then we must show that
$\underset{n\to \infty}{\lim}[\bo^{(n)}] = [\bo]$. Note that we suppress here our usual assumption $\om_r= 1$, which would
conflict with the normalization $g_s = \underset{n\to \infty}{\lim} g_s^{(n)} = 1$.
 \item[(iii)] If $s=r$, we are done. This follows from the fact that $j:\: \Ga \setminus \Om^r \stackrel{\cong}{\lra} M^r$ is an
isomorphism of analytic spaces, thus a strong homeomorphism. Hence we may suppose that $s <r$.
 \item[(iv)] Consider the Newton polygon ${\rm NP}(\phi)$ of $\phi_T$, i.e., the lower convex hull of the vertices 
 $(q^{i},v(g_i))$ with $0 \leq i \leq s$ in the plane (see \cite{30} II Sect. 6). If ${\rm rk}(\phi^{(n)})\leq s$ for 
 $n\gg 0$, then in fact ${\rm rk}(\phi^{(n)}) = s$ for $n \gg 0$, and we are ready as in (iii). Therefore, possibly restricting
 to a subsequence, we may assume that ${\rm rk}(\phi^{(n)}) >s$ for $n \gg 0$. Then if $\phi^{(n)}$ is sufficiently close
 to $\phi$, the Newton polygon ${\rm NP}(\phi^{(n)})$ agrees with ${\rm NP}(\phi)$ from the leftmost vertex $(1,-1)$ up to
 $(q^s,0)$ and, since  $g_{s+1}^{(n)}, \ldots,g_r^{(n)}$ tend to zero:\vspace{0.2cm}\\
(2.3.1) The slope of ${\rm NP}(\phi^{(n)})$ right to $(q^s,0)$ tends to infinity if $n \to \infty$.  
 \item[(v)] Considering \cite{18} (3.3), (3.4), (3.5), the assertion (2.3.1) implies that the quotient
 $|\om_{r-s}^{(n)}|/|\om_{r-s+1}^{(n)}|$ (i.e., the quotient of absolute values of the $(s+1)$-th divided by the $s$-the element 
 of our SMB $\{\om_r^{(n)},\om_{r-1}^{(n)},\ldots\}$) tends to infinity with $n\to \infty$.
 \item[(vi)] Let ${^s\phi}^{(n)}$ be the rank-$s$ Drinfeld module that corresponds to the lattice 
 ${^s\La}^{(n)} = A\om_r^{(n)} + \cdots+ A\om_{r-s+1}^{(n)}$, with 
  $${^s\phi}_T^{(n)}(X) = \sum_{0 \leq i \leq s} s_{g_i}^{(n)} X^{q^{i}},\: {^s\fg}^{(n)} := 
   ({^sg_1}^{(n)}: \ldots:{^sg_s}^{(n)}:0: \ldots:0).$$ 
Then $\underset{n\to \infty}{\lim} (g_i^{(n)} - {^sg}_i^{(n)}) = 0$, as follows from (v). (The analogous statement for the
Eisenstein series $E_{q^{i}-1}$, to wit
 $$\lim_{n\to \infty}(E_{q^{i}-1}(\La^{(n)})-E_{q^{i}-1}({^s\La}^{(n)})) = 0,$$
 is obvious; then we use the fact (2.1.6) that the $g_i$ are polynomials in the $E_k$.)
  \medskip
  
  \noindent
 As $g_i^{(n)} \lra g_i$ for $1 \leq i \leq r$, we find ${^s\fg}^{(n)} \lra \fg$ and therefore ${^s\phi}^{(n)} \lra \phi$ in
 $\overline{M}^s = {\rm Proj}\, C_{\infty}[g_1,\ldots,g_s]$ with respect to the strong topology.
 \item[(vii)] Denote by ${^s\bo}^{(n)}$ the point $(0:\ldots:0:\om_{r-s+1}^{(n)}:\ldots:\om_r^{(n)})$ in
 $\Om_{V_s} \hookrightarrow \overline{\Om}^r$. Applying (iii) with $r$ replaced by $s$ and using the 
 identification $\Om^s \stackrel{\cong}{\lra} \Om_{V_s}$,
  $$\lim_{n\to \infty} [{^s\bo}^{(n)}] = [\bo]\leqno{(2.3.2)}$$
holds in ${\rm GL}(s,A)\setminus \Om_{V_s} \hookrightarrow \Ga \setminus \overline{\Om}^r$, where $[~.~ ]$ is the
class modulo $\Ga$. But (2.3.2) together with (v) means that $[\bo^{(n)}]$ tends strongly to $[\bo]$ in 
$\Ga\setminus \overline{\Om}^r$.
   \end{itemize}
 \end{proof}  
  
(2.4) We want to give similar descriptions for the quotients $\Ga'\setminus \overline{\Om}^r$ and $\Ga'\setminus \overline{\Psi}^r$, 
where $\Ga' \subset \Ga$ is a congruence subgroup and $r \geq 2$. This turns out, however, to be much more difficult. We
restrict to deal with the case where $\Ga' = \Ga(N)$, the full congruence subgroup of level $N$. 
 \bigskip
 
Fix a monic $N\in A$ of degree $d \geq 1$, and write the $N$-th division polynomial of the Drinfeld module
$\phi^{\bo}$ ($\bo \in \Psi^r$) as
 $$\phi_N^{^\bo}(X) = \sum_{0 \leq i \leq rd} \ell_i(N,\bo)X^{q^{i}}$$
with $\ell_0(N,\bo) = N$, $\ell_{rd}(N,\bo) = \Delta(\bo)^{(q^{rd}-1)/(q^r-1)}$, where  $\Delta(\bo) = g_r(\bo)$ is the
discriminant function and, more generally, all the coefficient functions $\ell_i(N,.)$ lie in  ${\bf Mod} = C_{\infty}[g_1,\ldots,g_r]$.
It satisfies
 $$\phi_N^{\bo}(X) = \Delta(\bo)^{(q^{rd}-1)/(q^r-1)} \prod_{(\fu \in N^{-1}A/A)^r}(X-e_{\bo}(\fu \bo)).\leqno{(2.4.1)}$$
Here $\fu$ runs through a system of representatives of the finite $A$-module $(N^{-1}A/A)^r$ and $\fu \bo   = 
\underset{1 \leq i \leq r}{\sum} u_i\om_i$. That is, the
 $$d_{\fu}(\bo) := e_{\bo} (\fu\bo) \leqno{(2.4.2)}$$
are the $N$-division points of $\phi^{\bo}$. It is known (see \cite{23} Proposition 2.7 or \cite{13} 3.3.5) that for $\fu \not= 0$,
 $$d_{\fu}(\bo) = E_{\fu}(\bo)^{-1} \leqno{(2.4.3)}$$
with the partial Eisenstein series of weight 1 (see (1.12))
 $$E_{\fu} := E_{1,\fu}.\leqno{(2.4.4)}$$
 We draw the conclusions
  \bigskip
  
(2.4.5) $E_{\fu}$ never vanishes on $\Psi^r$ and $\Om^r$;
  \bigskip
  
(2.4.6) The coefficient $\ell_i(N,\bo)$ may be expressed as a homogeneous polynomial in the $E_{\fu}$ ($\fu \not= 0$);
more precisely, 
 $$\ell_i(N,\bo) = Ns_{q^{i}-1} (E_{\fu}(\bo)~|~0 \not= \fu \in (N^{-1}A/A)^r),$$
where $s_k$ is the $k$-th elementary symmetric polyomial. 
 \medskip
 
We let $\MT(N)$ be the index set
 $$\MT(N) := (N^{-1}A/A)^r \setminus \{0\}.\leqno{(2.4.7)}$$
(2.5) As the definition of fields of meromorphic functions on non-complete analytic spaces requires some boundary
conditions, we make the following {\it ad hoc} definitions. They are motivated from the fact that the analytic spaces
$\widetilde{M}^r$, $\widetilde{M}^r(N)$, $M^r$, $M^r(N)$ appearing below are actually $C_{\infty}$-varieties (see
(2.2.1), (2.2.5) and Remark 2.7) and that by GAGA the algebraic and the analytic function fields of projective 
$C_{\infty}$-varieties agree.
 \bigskip
 
(2.5.1) The function field of $\widetilde{M}^r = \Ga\setminus \Psi^r$ is 
  $$\widetilde{\MF}_r = \widetilde{\MF}(1) := C_{\infty}(g_1,\ldots,g_r);$$
(2.5.2) The function field of $\widetilde{M}^r(N):= \Ga(N) \setminus \Psi^r$ is $\widetilde{\MF}_r(N)$, the field
of those meromorphic functions on $\widetilde{M}^r(N)$ which are algebraic over $\widetilde{\MF}_r$;
 \bigskip
  
(2.5.3) The function field of $M^r = \Ga \setminus \Om^r$ is 
 $$\MF_r = \MF_r(1) := C_{\infty}(g_1,\ldots,g_r)_0,$$
the subfield of isobaric elements of weight 0 of $\widetilde{\MF}_r$;
 \bigskip  
  
(2.5.4) The function field of $M^r(N):= \Ga(N)\setminus \Om^r$ is $\MF_r(N)$, the field of meromorphic functions on
$M^r(N)$ algebraic over $\MF_r$. 
 \bigskip
 
{\bf 2.6 Proposition.}  
{\it 
\begin{itemize}
 \item[(i)] The field $\widetilde{\MF}_r(N)$ is generated over $C_{\infty}$ by the Eisenstein series $E_{\fu} = E_{1,\fu}$
 ($\fu \in \MT(N)$). It is galois over $\widetilde{\MF}_r$ with Galois group
  $$\widetilde{G}(N) := \{\ga \in {\rm GL}(r,A/M)~|~\det \ga \in \F^*\}.$$
   \item[(ii)] The field $\MF_r(N)$ is generated over $C_{\infty}$ by the functions 
   $E_{\fu}/E_{\fv}$ ($\fu,\fv  \in \MT(N)$). 
   It is galois
   over $\MF_r$ with group $G(N):= \widetilde{G}(N)/Z$. Here $Z \cong \F^*$ is the subgroup of $\widetilde{G}(N)$ of
   scalar matrices with entries in $\F^*$.
 \end{itemize}
 }
 
 \begin{proof}
  \begin{itemize}
  \item[(i)] As $\Ga$ acts without fixed points on $\Psi^r$, $\widetilde{M}^r(N) = \Ga(N)\setminus \Psi^r$ is an \'etale 
  Galois cover of $\widetilde{M}^r = \Ga \setminus \Psi^r$ with group $\Ga/\Ga(N) \stackrel{\cong}{\lra} \widetilde{G}(N)$.
  Now $E_{\fu}$ is $\Ga(N)$-invariant and, as (2.4.1) and (2.4.3) show, algebraic over $\widetilde{F}_r$, i.e., 
  $E_{\fu} \in \widetilde{\MF}_r(N)$. 
  Furthermore, the relation $E_{\fu}(\ga\bo) = E_{\fu \ga}(\bo)$ for $\bo \in \Psi^r$, $\ga \in \Ga$ implies $\ga = 1$ if
  $\ga\in  \widetilde{G}(N)$ fixes all the $E_{\fu}$. Therefore, $\widetilde{\MF}_r(N) = \widetilde{\MF}_r(E_{\fu}~|~\fu \in \MT(N))$
  by Galois theory. In view of (2.4.6) the coefficient functions $\ell_i(N,.)$ and therefore (by the well-known commutation 
  relations between the $g_i(~.~)$ and the $\ell_i(N,.)$) also the $g_i$ are polynomials in the $E_{\fu}$. Thus in fact
  $\widetilde{\MF}_r(N) = C_{\infty}(E_{\fu}~|~\fu \in \MT(N))$. 
   \item[(ii)] The argument for $\MF_r(N)$ is similar. The quotient $G(N) = \widetilde{G}(N)/Z$ by $Z$ as a Galois group
   comes from the fact that $Z$ acts trivially on $\Om^r$.
   \end{itemize}
  \end{proof}
  
 {\bf Remark.} By (2.4.3) we may also write $\widetilde{F}_r(N) = C_{\infty}(d_{\fu}~|~\fu \in \MT(N))$ and $\MF_r(N) = 
 C_{\infty}(d_{\fu}/d_{\fv}~|~\fu,\fv \in \MT(N))$.
  \bigskip
  
 {\bf 2.7 Remark.} As is well known, the smooth analytic space $M^r(N)=\Ga(N)\setminus \Om^r$ is strongly related with the
 moduli scheme $M^r(N)/K$ of Drinfeld $A$-modules of rank $r$ with a
 structure of level $N$ (\cite{11}, \cite{10}, \cite{15}). Let $K(N) \subset C_{\infty}$ be the field extension of $K$ generated by
 the $N$-division points of the Carlitz module. Then $K(N)/K$ is finite abelian with group $(A/N)^*$ and ramification 
 properties similar to those of cyclotomic extensions of $\mathbb Q$ \cite{26}. Let $K_+(N) \subset K(N)$ be the fixed field of 
 $\F^* \hookrightarrow (A/N)^*$, the ``maximal real subextension'' of $K(N)|K$. Then $K_+(N)$ is contained in 
 $\MK_r(N)=K(E_{\fu}/E_{\fv}~|~\fu,\fv \in \MT(N))$, and is actually the algebraic closure of $K$ in $\MK_r(N)$.
 Now $M^r(N)/K$ is a smooth $K$-scheme with function field $\MK_r(N)$, whose set of $C_{\infty}$-points (in fact, its 
 analytification over $C_{\infty}$) is given by
  $$(M^r(N)/K)(C_{\infty}) \stackrel{\cong}{\lra} \underset{\sigma}{\stackrel{\bullet}{\bigcup}}M^r(N)_{\sigma} =
 \underset{\sigma}{\stackrel{\bullet}{\bigcup}}  (\Ga(N)\setminus \Om^r)_{\sigma},$$ 
where $\sigma$ runs through the set of $K$-embeddings of $K_+(N)$ into $C_{\infty}$, i.e., the Galois group
${\rm Gal}(K_+(N)|K) = (A/N)^*/\F^*$. Correspondingly, the analytification of $M^r(N)/K \underset{K_+(N)}{\times} C_{\infty}$ is
$M^r(N)$, which justifies our notation $M^r(N)$ for $\Ga(N)\setminus \Om^r$.  In the language of pre-Grothendieck algebraic
geometry, $M^r(N)/K$ is ``defined over $K_+(N)$''.
  \bigskip

The group ${\rm GL}(r,A/N)$ acts naturally on the set of $N$-level structures of a fixed Drinfeld module of rank $r$, thus
on $M^r(N)/K$, which identifies $M^r/K = ({\rm Proj}\,K[g_1,\ldots,g_r])_{(g_r\not=0)}$ with the quotient of $M^r(N)/K$ by this
group. Moreover, the action is compatible with 
that of $G(N) = \Ga/\Ga(N)Z \hookrightarrow {\rm GL}(r,A/N)/Z$ on the components $M^r(N)_{\sigma}$. All of this may
be transferred to the spaces $\widetilde{M}^r(N) = \Ga(N)\setminus \Psi^r$ and their function fields $\widetilde{\MF}_r(N)$. 
As we don't really need it, we omit the details.
 \bigskip
 
In the sequel of the paper, we will construct a compactification $\overline{M}^r(N)$ of $M^r(N)$ (and, similarly, a 
horizontal compactification of $\widetilde{M}^r(N)$), i.e., a projective $C_{\infty}$-variety $\overline{M}^r(N)$ with
set of $C_{\infty}$-points $\overline{M}^r(N)(C_{\infty}) = \Ga(N)\setminus \overline{\Om}^r$, into which $M^r(N)$ embeds
as a dense open subvariety, and compatible with the above-described group actions. 
 \vspace{0.5cm}
 
{\bf 3. The boundary components}
 \vspace{0.3cm} 

From now on, we assume that $r \geq 2$.
 \bigskip

The set $\mathfrak U_s$ of $s$-dimensional subspaces $U$ of $V = K^r$ is in canonical bijection with
${\rm GL}(r,K)/P_s(K)$ through
 $$\begin{array}{rcl}
  {\rm GL}(r,K)/P_s(K) & \stackrel{\cong}{\lra} & \mathfrak U_s.\\
   \ga & \longmapsto & V_s\ga^{-1}
    \end{array}\leqno{(3.1)}$$

As the action of $\Ga$ on $\mathfrak U_s$ is transitive, we may replace the left hand side with 
$\Ga / \Ga\cap P_s(K)$. Let $(A/N)^r_{\rm prim}$ be the set of primitive elements of $(A/N)^r$, that is, of elements
that belong to a basis of the free $(A/N)$-module $(A/N)^r$. Then, as is easily verified, the map
 $$\Ga(N) \setminus \Ga/\Ga \cap P_{r-1}(K) \lra (A/N)^r_{\rm prim}/\F^* =: \MC_r(N)\leqno{(3.2)}$$
that associates with the double class of $\ga \in \Ga$ the first column of $\Ga$ (evaluated modulo $N$, and modulo
the scalar action of $\F^*$) is well-defined and bijective. Together with (3.1) we find that the space of orbits on
$\mathfrak U_{r-1}$ of $\Ga(N)$ is 
 $$\Ga(N)\setminus \mathfrak U_{r-1} \stackrel{\cong}{\lra} \MC_r(N).\leqno{(3.3)}$$
This allows us to describe the components of codimension 1 of 
 $$\Ga(N)\setminus \overline{\Psi}^r = \underset{1 \leq s \leq r}{\stackrel{\bullet}{\bigcup}} \:\:
 \underset{U\in \Ga(N)\setminus\mathfrak U_s}{\stackrel{\bullet}{\bigcup}} \Ga_U(N)\setminus \Psi_U \leqno{(3.4)}$$
 and, analogously, of $\Ga(N)\setminus \overline{\Om}^r$. Here $\Ga_U = \{\ga \in \Ga~|~U \ga = U\}$, which acts from the
 left on $\Psi_U$ (see (1.5)), and $\Ga_U(N) := \Ga_U \cap \Ga(N)$. We put
  $$\widetilde{M}_U(N) := \Ga_U(N)\setminus \Psi_U \mbox{ and } M_U(N) := \Ga_U(N)\setminus \Om_U \leqno{(3.5)}$$
and call the components $\widetilde{M}_U(N)$ resp. $M_U(N)$ with $\dim(U) = r-1$ the {\em cuspidal divisors} or simply
the {\em cusps} of $\Ga(N)\setminus \overline{\Psi}^r$ (or $\Ga(N) \setminus \overline{\Om}^r$, or of $\Ga(N)$). Each
of these sets is in canonical bijection with $\MC_r(N)$. For later use we specify a system of representatives, namely the
set  $S$ of monic elements of $(A/N)^r_{\rm prim}$. Here, some $\bn= (n_1,\ldots,n_r) \in (A/N)^r_{\rm prim}$ is {\em monic} if the 
first non-vanishing $n_i$ has a monic representative $n'_i \in A$ of degree less than $d = \deg\,N$.
 \bigskip
 
The cardinality $c_r(N)$ of $\MC_r(N)$ is an easy arithmetic function of $N$ and $r$, given by the following formula.
 \bigskip
 
{\bf 3.6 Lemma.} {\it Let $N = \underset{1\leq i \leq t}{\prod} \fp_i^{s_i}$ be the decomposition of $N$ into powers of 
different primes $\fp_i$
of $A$. Write $q_i=q^{\deg\,\fp_i}$. Then
 $$c_r(N) = (q-1)^{-1} \prod_{1 \leq i \leq t} (q_i^r-1) q_i^{(s_i-1)r}.$$}
  
 \begin{proof}
 We must determine $\#(A/N)^r_{\rm prim}$, which by the Chinese Remainder Theorem is multiplicative. So we may
 assume that $t=1$, $N = \fp_1^{s_1}$ with some prime $\fp_1$, $q_1 = q^{\deg\,\fp_1}$. Then 
  $$\#(A/N)^r_{\rm prim} = \#(A/\fp_1)_{\rm prim}^r \#(\fp_1/\fp_1^{s_1})^r,$$
  as some element of $(A/N)^r$ is primitive if and only if its reduction $\bmod\,\fp_1$ is. Now $\#(A/\fp_1)^r_{\rm prim} =
  q_1^r-1$ and $\#(\fp_1/\fp_1^{s_1})^r = q_1^{(s_1-1)r}$, and we are done.
 \end{proof} 
  \bigskip
  
In the case of smaller dimension $s < r-1$ we get a similar description of $\Ga(N)\setminus \mathfrak U_s$, which is in
$1-1$-correspondence with the set of $(r-s)$-subsets of $(A/N)^r$ that are part of an $(A/N)$-basis of $(A/N)^r$, modulo the
action of the group $\{\ga \in {\rm GL}(r-s,A/N)~|~\det\,\ga \in \F^*\}$. We leave the details to the reader, as we will only
need the case $s=1$. Here, likewise,
 $$\Ga(N)\setminus \mathfrak U_1 \stackrel{\cong}{\lra} \Ga(N)\setminus \Ga/\Ga\cap P_1(K) \stackrel{\cong}{\lra} \MC_r(N),
  \leqno{(3.7)}$$
 where the double class of $\ga \in \Ga$ is mapped to the last row vector of $\ga$ (reduced modulo $N$, and modulo the 
 action of $\F^*$). In particular, 
  $$\#(\Ga(N)\setminus \mathfrak U_1) = \#\MC_r(N) = c_r(N),$$
which of course could also be seen via the duality of projective spaces over the finite ring $A/N$.
 \vspace{0.5cm}
 
{\bf 4. Behavior of Eisenstein series at the boundary.}
 \vspace{0.3cm} 
 
In the whole section, $N$ is a fixed monic element of $A$ of degree $d \geq 1$. 
 \bigskip
 
(4.1) Let $\fu$ be an element of $\MT(N) = (N^{-1}A/A)^r \setminus \{0\}$. We start with the relation from (1.12)
 $$E_{k,\fu}(\ga(U,i)) = E_{k,\fu\ga}(U,i)\leqno{(4.1.1)}$$
for $\ga \in \Ga$, $(U,i) \in \overline{\Psi}^r$. Suppose that $U = V_s \cdot \ga^{-1}$ with some $1 \leq s < r$ and
$\ga \in \Ga$. Now we read off from (4.1.1):
 \bigskip
 
(4.1.2) The vanishing behavior of $E_{k,\fu}$ around the boundary component $\Psi_U$ is the same as the behavior
of $E_{k,\fu\ga}$ around the standard component $\Psi_{V_s}$.
 \bigskip
 
We say that
 \bigskip
 
(4.1.3) $\fu$ {\em belongs to} $U$ if $\fu \in \MT(N) \subset K^r/A^r = V/L_V$ is represented by an element of $U \subset V$.
 \bigskip
 
In view of $\fu\ga = \fu$ for $\ga \in \Ga(N)$, this property depends only on the $\Ga(N)$-orbit of $U$.
 \bigskip
 
{\bf 4.2 Proposition.} {\it 
 \begin{itemize}
  \item[(i)] Suppose  that $\fu$ does not belong to $U$. Then $E_{k,\fu}$ vanishes identically on $\Psi_{U}$.
  \item[(ii)] If $\fu$ belongs to $U$, then $E_{k,\fu}$ restricts to $\Psi_U \cong \Psi^s$ like an Eisenstein series $E_{k,\fu'}$
of rank $s$, $\fu' \in \MT_s(N) = (N^{-1}A/A)^s \setminus \{0\}$. In particular, it doesn't vanish identically on $\Psi_U$. 
  \end{itemize}     
  } 
  \begin{proof}
In view of (4.1.2), it suffices to verify the assertions for $U=V_s$. Suppose that $\fu = (u_1,\ldots,u_r)$ does not belong
to $V_s$, that is, $u_i \not= 0$ for some $i$ with $1 \leq i \leq r-s$. Let $\bo = (\om_1,\ldots,\om_r)$ be an element of the
fundamental domain $\widetilde{\BF}$ described in (1.15). We have 
 $$E_{k,\fu}(\bo) = N^{-k} \sum_{\fa \in A^r \atop \fa \equiv N\fu (\bmod N)} (a_1 \om_1+ \cdots + a_r\om_r)^{-k} = 
 N^{-k} \sum_{\fa} (\fa\bo)^{-k}.$$
In each term, $a_i \not= 0$, which by (1.8.1) forces that $(\fa\bo)^{-k}$ tends to zero, uniformly in the $\fa$, if $\bo$
approaches $\Psi_{V_s}$. That is, $E_{k,\fu}(\bo) \lra 0$, and $E_{k,\fu} \equiv 0$ on $\Psi_{V_s}$. Suppose that
$\fu$ belongs to $V_s$. As before, each term $(\fa\bo)^{-k}$ tends to zero uniformly in $\fa$, as long as at least one
of $a_1,\ldots,a_{r-s} \not=0$. Therefore, $\lim\,E_{k,\fu}(\bo) = E_{k,\fu'}(\bo')$ with $\fu' = (u_{r-s+1},\ldots ,u_r)$ if
$\bo$ tends to $(0,\ldots,0,\om'_{r-s+1},\ldots,\om'_r) = (0,\ldots,0,\bo')$. 
   \end{proof}  
 We define the space
  $${\rm Eis}_k(N) := \sum_{\fu \in (N^{-1}A/A)^r} C_{\infty} E_{k,\fu}\leqno{(4.3)}$$
of Eisenstein series of weight $k$ and level $N$.
 \bigskip
 
{\bf 4.4 Lemma.} {\it    
The vector space ${\rm Eis}_k(N)$ is generated by the $E_{k,\fu}$ with $\fu$ primitive of level $N$ (i.e., $N'\fu \not= 0$
for each proper divisor $N'$ of $N$) and even by $E_{k,\fu}$ ($\fu \in N^{-1}S$), where $S$ is the set of representatives
for $\MC_r(N) = (A/N)^r_{\rm prim}/\F^*$ given in $(3.5)$.} 
 
  \begin{proof} Let $N'$ be a monic divisor of $N$, where $N' = 1$ is allowed. Then the distribution relation
   $$(N/N')^k \sum_{(N/N')\fu = \fv} E_{k,\fu} = E_{k,\fv}\leqno{(4.4.1)}$$
 holds, where $\fv \in (N'^{-1}A/A)^r$ and $E_{k,0} = E_k$ is the Eisenstein series without level. It shows that Eisenstein
 series of lower level $N'$ may be expressed as linear combinations of those with level $N$. The result now follows from
  $$E_{k,c\fu} = c^{-k}E_{k,\fu} \quad (c \in \F^*), \leqno{(4.4.2)}$$
 which is immediate from the definition of $E_{k,\fu}$. 
  \end{proof}  
  
 (4.5) We will show that these are all the relations between Eisenstein series of weight $k$, following the strategy of Hecke
 in \cite{27}, which has been introduced to the function field situation in the case of $r=2$ by Cornelissen \cite{9}. Let
 \bigskip
 
 $$F_{k,\fu}(\bo) := N^{-k} \sum_{\fa \in A^r\,{\rm primitive} \atop \fa \equiv N {\fu}\,(\bmod N)} (\fa \bo)^{-k}
  \leqno{(4.5.1)}$$
 be the partial sum of $E_{k,\fu}(\bo)$ with primitive $\fa$, i.e., where $\fa = (a_1,\ldots,a_r)$ satisfies 
 $\underset{1 \leq i \leq r}{\sum} A a_i = A$. Then:
  \bigskip
  
(4.5.2) The {\em restricted Eisenstein series} $F_{k,\fu}$ is well-defined as a function on $\Psi^r$ of weight $k$ and invariant
under $\Ga(N)$. Like $E_{k,\fu}$, it satisfies the functional equation
 $$F_{k,\fu}(\ga\bo) = F_{k,\fu\ga}(\bo)$$
 under $\ga \in \Ga$.
  \bigskip
  
Let $\mu:\,A \lra  \{0,\pm 1\}$ be the M\"obius function:\\
 $\mu(a) = (-1)^n \mbox{ if } a = \epsilon \underset{1 \leq j \leq n}{\prod} \fp_j$ with $n$ different monic primes $\fp_j$ of $A$
 and $\epsilon \in \F^*$, and zero otherwise. (As the empty product evaluates to 1, $\mu(a) = 1$ if $a \in \F^*$.) Then
$\underset{b\,{\rm monic},\, b|a}{\sum} \mu(b) = 1$ if $a \in \F^*$ and 0 otherwise, and the usual formalism holds. M\"obius inversion
yields
 $$F_{k,\fu}(\bo) = \sum_{t\in(A/N)^*} \sum_{a\in A\,{\rm monic}\atop
 at \equiv 1(\bmod N)} \mu(a)a^{-k} E_{k,t\fu}(\bo). \leqno{(4.5.3)}$$
In particular, $F_{k,\fu}$ lies in ${\rm Eis}_k(N)$, and so has a strongly continuous extension to $\overline{\Psi}^r$. We deduce
that 
 $$\dim({\rm Eis}_k(N)) \geq \dim(\sum_{\bn \in \MT(N)} C_{\infty} F_{k,\fu}). \leqno{(4.5.4)}$$
Recall that by (3.7) the set of 1-dimensional boundary components of $\Ga(N)\setminus \overline{\Psi}^r$ is in $1-1$-correspondence
with $\MC_r(N)$, or with its set $S$ of representatives in (3.5). We let $\Psi_{\bn}^1 \cong \Psi^1$ be the component
corresponding to $\bn \in S$.
 \bigskip
 
{\bf 4.6 Proposition.} {\it
 \begin{itemize}
  \item[(i)] Given $\bn \in S$ there exists a unique $\bn' \in S$ such that $F_{k,\bn'/N}$ doesn't vanish at $\Psi_{\bn}^1$.
   \item[(ii)] The rule $\bn \lra \bn'$ establishes a permutation of $S$.
  \end{itemize} 
  }
  \begin{proof}
First we note that $\Psi_{\bn}^1$, where $\bn = (0,\ldots,0,1)$, equals $\Psi_{V_1}$. In view of (4.5.2) and the transitivity
of $\Ga$ on the set $\{\Psi_{\bn}^1~|~\bn \in S\}$, it is enough to show that there is a unique $\bn'$ such that
$F_{k,\bn'/N}$ doesn't vanish at $\Psi_{V_1}$. 
 \bigskip
 
Consider a term $(\fa\bo)^{-k} = (a_1\om_1+ \cdots+ a_r\bo_r)^{-k}$ of $N^kF_{k,\bn'/N}$, where $\gcd(a_1,\ldots,a_r) = 1$ and
$\bo \in \tilde{\BF}$ as in the proof of (4.2). If one of $a_1,\ldots,a_{r-1}$ doesn't vanish then  $\lim(\fa\bo)^{-k} = 0$, uniformly in
$\fa$, if $\bo$ tends to $(0,\ldots,0,\om_r)$. Hence
 $$\lim F_{k,\bn'/N}(\bo) = \lim N^{-k} \underset{\fa \equiv \bn'(\bmod N)}{\sum_{\fa\in A^r {\rm \,primitive}\atop a_1\ldots a_{r-1}=0}}
 (\fa\bo)^{-k},$$
 which can be non-zero only if $a_r \in \F^*$ (and in fact $a_r = 1$, as $\bn'$ is required to be monic). This implies that
 $\bn' = (0,\ldots,0,1)$. Conversely, that choice of $\bn'$ gives $F_{k,\bn'/N}(0,\ldots,0,1) = N^{-k}\om_r^{-k} \not= 0$. That is,
$\bn' = \bn = (0,\ldots,0,1)$ is as wanted.
  \end{proof}
  
 {\bf 4.7 Corollary.} {\it The restricted Eisenstein series $F_{k,\fu}$, where $\fu \in N^{-1}S$, are linearly independent and form a
 basis of ${\rm Eis}_k(N)$. The dimension of ${\rm Eis}_k(N)$ equals $\#(S) = c_r(N)$.}

  \begin{proof}
  (4.4) + (4.5.4) + (4.6)
   \end{proof}  
   \medskip
   
{\bf 4.8 Proposition.} {\it
  \begin{itemize}
   \item[(i)] The Eisenstein series $E_{\fu} = E_{1,\fu}$ ($u \in \MT(N)$) of weight $1$ and level $N$ separate points of
$\Ga(N)\setminus \overline{\Psi}^r$. That is, if $\bo, \bo' \in \overline{\Psi}^r$ satisfy $E_{\fu}(\bo) = E_{\fu}(\bo')$ for all
$\fu \in \MT(N)$, then there exists $\Ga \in \Ga(N)$ such that $\bo' = \ga \bo$. 
 \item[(ii)] The same statement for $\Ga(N)\setminus \overline{\Psi}^r$ replaced with $\Ga(N)\setminus \overline{\Om}^r$.
 \end{itemize} 
   } 
   
    \begin{proof}
We start with the identity
  $$NX \prod_{\fu\in \MT(N)}(1-E_{\fu}X) = \phi_N(X) = NX + \sum_{1 \leq i \leq rd} \ell_i(N)X^{q^{i}}, \leqno{(4.8.1)}$$
which comes from (2.4.1) and (2.4.3). Here the right hand side is the $N$-division polynomial of the general Drinfeld module
$\phi$ of rank $\leq r$, which lives above $\overline{\Psi}^r$. The coefficients $\ell_i(N)$, where $1 \leq i \leq rd$ ($d = \deg N$), are 
$\Ga$-invariant functions on $\overline{\Psi}^r$ of weights $q^{i}-1$. 
\bigskip

Hence the data $\{E_{\fu}(\bo)~|~\fu \in \MT(N)\}$ determines the values on $\bo$ of the $\ell_i(N,\bo)$ and therefore, 
taking the known relations between the $\ell_i(N)$ and the $g_i$ ($1 \leq j \leq r$) into account, the coefficients 
$g_1(\bo),\ldots,g_r(\bo)$ of the $T$-division polynomial 
 $$\phi_T^{\bo}(X) = TX+\sum_{1 \leq j \leq r} g_j(\bo)X^{q^{j}}$$
of the Drinfeld module $\phi^{\bo}$ that corresponds to $\bo \in \overline{\Psi}^r$. Suppose that $\phi^{\bo}$ has rank
$s$ ($1 \leq s \leq r$), that is, $g_s(\bo) \not= 0$, $g_{s+1}(\bo) = \cdots = g_r(\bo) = 0$.
 \bigskip
 
If \fbox{$s=r$} then $\phi^{\bo}$ determines an $A$-lattice in $C_{\infty}$ of rank $r$, hence a point $\bo \in \Psi^r$ up to
the action of $\Ga$. That is, if $E_{\fu}(\bo) = E_{\fu}(\bo')$ for all $\fu$ then $\bo' = \ga \bo$ with  some $\ga \in \Ga$. The 
relation $E_{\fu}(\ga\bo) = E_{\fu\ga}(\bo)$ moreover shows that $\fu\ga = \fu$ for all $\fu$, that is, $\ga$ lies in fact in 
$\Ga(N)$.
 \bigskip
 
If \fbox{$s < r$}, the $A$-lattice corresponding to $\phi^{\bo}$ has rank $s$ and is given by an embedding $i:\: U \hookrightarrow
C_{\infty}$ of some $s$-dimensional subspace $U$ of $K^r$ (i.e., of $L_U = A^r \cap U \hookrightarrow C_{\infty}$). Now
Proposition 4.2 allows to determine $U$ up to $\Ga(N)$-equivalence. Choose one such $U$; then $\bo \in \Psi_U$ is determined
through $\phi^{\bo}$ and the $E_{\fu}(\bo)$ up to an element of $\Ga_U$ and, in fact (with the same argument as in the case
\fbox{$s=r$}), up to an element of $\Ga_U(N) = \Ga_U \cap \Ga(N)$. This shows (i); the proof of (ii) is identical.
     \end{proof}
     \vspace{0.5cm}
     
{\bf 5. The projective embedding.}
 \vspace{0.3cm}
 
In this section, we show that $\Ga(N)\setminus \overline{\Om}^r$ is the set of $C_{\infty}$-points of a closed subvariety
of some projective space. This allows us to endow it with the structure of projecive variety, which then will be labelled
with the symbol $\overline{M}^r(N)$, the {\em Eisenstein compactification} of $M^r(N)$.
 \bigskip
 
Throughout, $N\in A$ of degree $d \geq 1$ is fixed.
 \bigskip
 
(5.1) We define the {\em Eisenstein ring of level} $N$, ${\bf Eis}(N)$, as the $C_{\infty}$-subalgebra of $\widetilde{\MF}_r(N)$
generated by the Eisenstein series $E_{\fu} = E_{1,\fu}$ of level $N$ and weight 1. It is graded with respect to weight:
its $k$-th piece ${\bf Eis}_k(N)$ is the $C_{\infty}$-space generated by monomials of degree $k$ in the $E_{\fu}$. In particular,
${\bf Eis}(N)$ is generated as an algebra by ${\bf Eis}_1(N) = {\rm Eis}_1(N)$, a vector space of dimension $c_r(N)$ 
(see (4.7)). We also let ${\bf Mod} = {\bf Mod}(1) = C_{\infty}[g_1,\ldots,g_r]$ be the graded algebra of modular forms of
type zero for $\Ga$ (\cite{20} 1.7).      
 \bigskip
 
{\bf 5.2 Proposition.}
{\it
 \begin{itemize}
  \item[(i)] ${\bf Eis}(N)$ contains the algebra ${\bf Mod}$;
   \item[(ii)] ${\bf Eis}(N)$ is integral over ${\bf Mod}$;
    \item[(iii)] ${\bf Eis}(N)$ contains all the Eisenstein series $E_{k,\fu}$ of  arbitrary weight $k$.
  \end{itemize} 
 }  
 
  \begin{proof}
   \begin{itemize}
\item[(i)] The argument in the proof of Proposition 4.8 shows that $g_1,\ldots,g_r \in {\bf Eis}(N)$.\\
\item[(ii)]The $E_{\fu}$ are the zeroes of the monic polynomial $N^{-1}X^{q^{rd}} \phi_N(X^{-1})$ with coefficients 
$N^{-1} \ell_i(N) \in {\bf Mod}$, where $\phi_N(X)$ is as in (4.8.1).\\
\item[(iii)] Let $\La$ be any rank-$r$ $A$-lattice in $C_{\infty}$ and $G_{k,\La}(X)$ be its $k$-th Goss polynomial 
(\cite{23} 2.17, \cite{16} 3.4). It is of shape
 $$G_{k,\La}(X) = \sum_{0 \leq i \leq k} a_i(\La)X^{k-i}, \leqno{(5.2.1)}$$
where $a_i$ is a modular form of weight $i$ and type $0$, that is, $a_i \in {\bf Mod}$. (This follows from \cite{16} 3.4(ii).)
The characteristic property of Goss polynomials ({\it loc. cit} 3.4(i)) implies
 $$E_{k,\fu} = G_{k,\La}(E_{1,\fu}).\leqno{(5.2.2)}$$
Now (iii) is a consequence of (5.2.1) and (5.2.2).  
 \end{itemize}
   \end{proof}
    
(5.3) We define $\mathbb P = \mathbb P(N)$ as the projective space $\mathbb P({\rm Eis}_1(N)^{\wedge})$    
associated with the dual vector space ${\rm Eis}_1(N)^{\wedge}$ of ${\rm Eis}_1(N)$. As a scheme,
 $$\mathbb P = {\rm Proj}\, R,\leqno{(5.3.1)}$$
where $R := {\rm Sym}({\rm Eis}_1(N))$ is the symmetric algebra on ${\rm Eis}_1(N)$. Consider the map
 $$j_N:\: \overline{\Om}^r \lra \mathbb P \leqno{(5.3.2)}$$
to the $C_{\infty}$-valued points of $\mathbb P$ that with the class of $(U,i)$ associates the class (up to scalars) of
the linear form $E_{\fu} \longmapsto E_{\fu}(U,i)$. Then:
 \bigskip
 
(5.3.3) $j_N$ is well-defined, as the $E_{\fu}$ have weight 1 and for each $(U,i)$ there exists $\fu$ such that
$E_{\fu}(U,i) \not= 0$.
 \bigskip
 
(5.3.4) $j_N(\ga(U,i)) = j_N(U,i)$ for $\ga \in \Ga(N)$, as the $E_{\fu}$ are $\Ga(N)$-invariant. Here and in the following,
we write $j_N(U,i)$ for $j_N$ (class of $(U,i)$).
 \bigskip
 
(5.3.5) As a map from $\Ga(N)\setminus \overline{\Om}^r$ to $\mathbb P$, $j_N$ is injective, due to Proposition 4.8.
 \bigskip
 
(5.4) The graded algebra $R$ is supplied with a canonical homomorphism
  $$\epsilon:\: R \lra {\bf Eis}(N),\leqno{(5.4.1)}$$  
 which is the identity on $R_1 = {\rm Eis}_1(N) = {\bf Eis}_1(N)$ and surjective, since ${\bf Eis}_1(N)$ generates 
 ${\bf Eis}(N)$. Let $J$ be the kernel of $\epsilon$. Since ${\bf Eis}(N)$ is a domain, $J$ is a (homogeneous) prime ideal
 of $R$, and in particular, saturated (\cite{25} p. 125). Then $j_N$ maps $\Ga(N)\setminus \overline{\Om}^r$ to the vanishing
 variety $V(J) \subset \mathbb P$ of $J$.
  \bigskip
  
(5.5) Let $x \in V(J)$ be given. The proof of Proposition 4.8 shows that there exists an element $(U,i)$ of $\overline{\Psi}^r$,
well-defined up to the action of $\Ga(N)$, such that $j_N(U,i) = x$. Viz., for  simplicity choose a representative
$\widetilde{x} \in {\rm Eis}_1(N)^{\wedge} = {\rm Hom}_{C_{\infty}}({\rm Eis}_1(N),C_{\infty})$ and put
$\widetilde{x}_{\fu} := \widetilde{x}(E_{\fu})$. Interpreting $\widetilde{x}_{\fu}$ as a value of $E_{\fu}$, the
$\widetilde{x}_{\fu}$ determine the values of the coefficient forms $g_1,\ldots,g_r$ as in (4.8), thus (if $g_r\not=0$)
a point $\widetilde{\bo} \in \Psi^r$ up to the action of $\Ga(N)$. The corresponding point $\bo \in \Om^r$ is independent of
the choice of $\widetilde{x}$ and serves the purpose. If $g_s \not= 0$, $g_{s+1} = \ldots = g_r = 0$ then, as in (4.8),
the $s$-dimensional $K$-space $U$ and its boundary coponent $\Psi_U$ is determined up to $\Ga(N)$-equivalence
by the (non-) vanishing of the $\widetilde{x}_{\fu}$. Choosing one such $U$, there exists an embeddding $i:\: U \hookrightarrow
C_{\infty}$, unique up to $\Ga_U(N)$, that fits the given data. Then the class of $(U,i)$ in $\overline{\Om}^r$ is as wanted.
That is, 
 $$j_N:\: \Ga(N)\setminus \overline{\Om}^r \stackrel{\cong}{\lra} V(J) \leqno{(5.6)}$$
is in fact bijective. Furthermore, the restriction of $j_N$ to a stratum $\Ga_U(N)\setminus \Om_U$ of $\Ga(N)\setminus 
\overline{\Om}^r$ is analytic with respect to the analytification of $V(J)$, as the $E_{\fu}$ are. 
 \bigskip
 
(5.7) Next, we consider the canonical morphism $\kappa:\: V(J) \lra \overline{M}^r$ defined as follows: Choose elements
$G_i\in R = {\rm Sym}({\rm Eis}_1(N))$ such that $\epsilon(G_i) = g_i$ ($1 \leq i \leq r$, see (5.4)). For given
$x \in V(J)$, $G_i$ may be evaluated on $\widetilde{x}$ (notation as in (5.5)), and we put 
 $$\kappa(x) = (G_1(\widetilde{x}): \ldots : G_r(\widetilde{x})),$$ 
which is independent of the choice of $\widetilde{x}$ above $x$, and of the choices of the $G_i$. Furthermore, the 
diagram
 $$\begin{array}{ccl}
  \Ga(N)\setminus \overline{\Om}^r & \stackrel{j_N}{\lra} & V(J) \vspace{0.2cm}\\
   \downarrow \pi & & \hspace*{0.3cm}\downarrow \kappa\vspace{0.2cm}\\
  \Ga\setminus \overline{\Om}^r & \stackrel{j}{\lra} & \overline{M}^r = {\rm Proj}\, C_{\infty}[g_1,\ldots,g_r]
   \end{array} \leqno{(5.7.1)}$$ 
commutes, where the left vertical arrow is the canonical projection $\pi$.
 \bigskip
 
{\bf 5.8 Proposition.} {\it $j_N$ is a strong homeomorphism. }
 \bigskip
 
\begin{proof}
This follows essentially from Theorem 2.3, that is, from the corresponding property of $j$. As in the proof of (2.3), 
$j_N$ is strongly continuous, so we must show that it is also an open map. Consider diagram (5.7.1), where $\pi$ and
therefore $j \circ \pi = \kappa \circ j_N$ are open. As $\kappa:\: V(J) \lra \overline{M}^r$ is set-theoretically the quotient map of the 
finite group $G(N) = \Ga/\Ga(N)\cdot Z$, which acts through homeomorphisms on $V(J)$, the openness of $\kappa \circ j_N$
implies the openness of $j_N$.
 \end{proof} 
By (5.6) and (5.8), we may use $j_N$ to endow $\Ga(N)\setminus \overline{\Om}^r$ with the structure of (the set of
$C_{\infty}$-points of) the projective subvariety $V(J)$ of $\mathbb P$, compatible with the analytic structures and the
strong topologies on both sides. By construction, $V(J)$ equals the projective variety associated with the graded algebra
${\bf Eis}(N)$. We collect what has been shown.
 \bigskip
  
{\bf 5.9 Theorem.} {\it
Let $N$ be a non-constant monic element of $A$.
 \begin{itemize}
  \item[(i)] The set $\Ga(N)\setminus \overline{\Om}^r$ is the set of $C_{\infty}$-points of an irreducible projective variety
  $\overline{M}^r(N)$ over $C_{\infty}$, the {\em Eisenstein compactification} of $M^r(N)$, which may be described as
  the variety ${\rm Proj}\,{\bf Eis}(N)$ associated with the Eisenstein ring ${\bf Eis}(N)$. It is a closed subvariety of the
  projective space $\mathbb P = \mathbb P({\rm Eis}_1(N)^{\wedge})$ attached to the dual of the vector space 
  ${\rm Eis}_1(N)$ of Eisenstein series of level $N$ and weight $1$, which has dimension $c_r(N)$. The open subvariety
  $M^r(N) = \Ga(N)\setminus \Om^r$ of $\overline{M}^r(N)$ is characterized as 
  $\{x \in \overline{M}^r(N)~|~E_{\fu}(x) \not= 0 \quad \forall \, \fu \in \MT(N)\}$.
   \item[(ii)] The set $\Ga(N)\setminus \overline{\Psi}^r$ is the set of $C_{\infty}$-points of an irreducible variety 
  $\overline{\widetilde{M}}^r(N)$ over $C_{\infty}$, which may be described as the variety ${\rm Spec}\,{\bf Eis}(N) \setminus \{I\}$,
  where $I$ is the irrelevant ideal of the graded ring ${\bf Eis}(N)$. It is a subvariety of the affine space attached to
  ${\rm Eis}_1(N)^{\wedge}$, endowed with an action of the multiplicative group $\mathbb G_m$, and such that
   $$\mathbb G_m \setminus \overline{\widetilde{M}}^r(N) \stackrel{\cong}{\lra} \overline{M}(N).$$
  \end{itemize}
 }
 
  \begin{proof}
 (i) has been shown above (see Proposition 4.2 for the last assertion), and the proof of (ii) is - mutatis mutandis -
 identical.
   \end{proof}
   
 {\bf 5.10 Remark.} We point out the following functorial properties of the construction of $\overline{M}(N)$ and
 $\overline{\widetilde{M}}^r(N)$.

  \begin{itemize}
   \item[(i)] It is compatible with level changes, to wit: Let $N'$ be a multiple of $N$ and $G(N,N')$ the quotient group
   $\Ga(N)/\Ga(N')$. The action of $\Ga(N)$ on $\Psi^r$ induces an action of $G(N,N')$ on $\Ga(N')\setminus \Psi^r = 
   \widetilde{M}^r(N')$ such that $G(N,N') \setminus \widetilde{M}^r(N') = \widetilde{M}^r(N)$. Further, the fixed
   space of $\Ga(N)$ in ${\rm Eis}_1(N')$ is ${\rm Eis}_1(N)$; hence $G(N,N')$ acts effectively on ${\rm Eis}_1(N')$ with
   fixed space
   ${\rm Eis}_1(N)$. Let $\widetilde{j}_N:\: 
   \Ga(N)\setminus \overline{\Psi}^r \hookrightarrow {\rm Eis}_1(N)^{\wedge}$ be the morphism analogous to $j_N$ and
   implicitly referred to in Theorem 5.9(ii). Then the diagram
  $$\begin{array}{ccc}
  \Ga(N')\setminus \overline{\Psi}^r & \stackrel{\widetilde{j}_{N'}}{\lra} & {\rm Eis}_1(N')^{\wedge}\vspace{0.2cm}\\
   \downarrow & & \hspace*{0.3cm}\downarrow \vspace{0.2cm}\\
  \Ga(N)\setminus \overline{\Psi}^r & \stackrel{\widetilde{j}_N}{\lra} & {\rm Eis}_1(N)^{\wedge}
   \end{array} \leqno{(5.10.1)}$$ 
is commutative and compatible with the action of $G(N,N')$, where the vertical arrows are the canonical projections.
In particular, the action of $G(N,N')$ on $\widetilde{M}^r(N')$ with quotient $\widetilde{M}^r(N)$ extends to
$\overline{\widetilde{M}}^r(N')$ with quotient
$\overline{\widetilde{M}}^r(N)$. Factoring out
the multiplicative group $\mathbb G_m$, we find similarly that $G(N,N')$ acts on $\overline{M}^r(N')$ with quotient
$\overline{M}^r(N)$.
 \item[(ii)] The construction of the Eisenstein compactification $\overline{M}^r(N)$ (and likewise of
 $\overline{\widetilde{M}}^r(N)$) is hereditary in the following sense. Let $\Om_U$ be a boundary component 
 ($U \in \mathfrak U_s$, $s < r$) and \\
  $$M_U(N):= \Ga_U(N)\setminus \Om_U \stackrel{\cong}{\lra} \{\ga \in {\rm GL}(s,A)~|~\ga \equiv 1 (\bmod N)\}
   \setminus \Om^s$$
its image in $\overline{M}^r(N)$. Then the Zariski closure $\overline{M}_U(N)$ of $M_U(N)$ in $\overline{M}^r(N)$ is composed of the $M_{U'}(N)$, where $U' \in \mathfrak U$ and $U' \subset U$, and is isomorphic with the variety $\overline{M}^s(N)$. This is
seen by assuming, without restriction, that $U = V_s$, in which case the description of $\overline{M}_{V_s}(N)$ is
identical with that of $\overline{M}^s(N)$.   
   \end{itemize}
  \bigskip
  
{\bf 5.11 Remark.} The idea of using Eisenstein series for a projective embedding of $M^r(N)$ is taken from \cite{28}.     
However, Kapranov's construction has the drawback that it fails to be canonical (it depends on the choice of a certain bound
$m_0$, see \cite{28} Proposition 1.12). Instead, our Proposition 4.8 assures that it suffices to consider Eisenstein series
of weight 1, which culminates in the canonical description $\overline{M}^r(N) = {\rm Proj}\,{\bf Eis}(N) \hookrightarrow 
\mathbb P({\rm Eis}_1(N)^{\wedge})$ with its functorial properties.
 \vspace{0.5cm}
 
{\bf 6. Tubular neighborhood of cuspidal divisors.}
 \vspace{0.3cm}
 
In this section we show that each point $x$ on a cuspidal divisor, i.e., on a boundary component $M_U(N)$ of
$\overline{M}^r(N)$ of codimension 1, possesses a neighborhood $Z$ isomorphic with $B \times W$, where $W$ is an open
admissible affinoid neighborhood of $x$ on $M_U(N)$ and $B$ a ball, and such that the map $\pi:\: Z \lra W$ derived
from the canonical projection $\pi_U:\: \overline{M}^r(N) \lra M_U(N)$ is the projection to the second factor.
 \bigskip

(6.1) As usual, it suffices to treat the case where $x$ is represented by $\bo^{(0)} \in \Om_{V_{r-1}} \stackrel{\cong}{\lra}
\Om^{r-1}$. For simplicity, we use the canonical isomorphism as an identification. We may further assume that
$\bo^{(0)}$ belongs to the fundamental domain $\BF'$ of $\Ga' = {\rm GL}(r-1,A)$ in $\Om^{r-1}$, that is,
$\bo^{(0)} = (0:\om_2^{(0)}: \ldots: \om_r^{(0)})$, where $\{1=\om_r^{(0)},\ldots,\om_2^{(0)}\}$ is an SMB of its lattice. Let
$X \subset \Om^{r-1}$ be the subspace
 $$X=\{\bo' = (\om'_2:\ldots:\om'_r = 1)~|~|\om'_i| = |\om_i^{(0)}|,\: 2 \leq i \leq r\}.\leqno{(6.1.1)}$$
Then, in fact, $X \subset \BF'$ and $X$ is an admissible open affinoid subspace, whose structure has been investigated
in \cite{20}, Theorem 2.4. (All of this collapses for $r=2$ to $X = \BF_1 = \Om^1 = \{{\rm point}\}$.) 
 \bigskip
 
We next put
 \medskip
 
 (6.1.2) $Y_c =$\\
$\{\bo \in \Om^r~|~\bo=(\om_1:\ldots:\om_r)~|~(\om_2:\ldots:\om_r) \in X~|~d(\om_1, \langle \om_2,\ldots,\om_r\rangle_{K_{\infty}})
 \geq c\}$
for some large $c$ in the value group $q^{\BQ}$ of $C_{\infty}$. Here $d(\om,\langle~.~\rangle_{K_{\infty}})$ is the
distance function to the $K_{\infty}$-space generated by $\om_2,\ldots,\om_r=1$. It is an admissible open subspace
of $\Om^r$. Note that $\bo \in Y_c$ in particular implies $|\om_1| \geq c$.  
 \bigskip
 
{\bf 6.2 Lemma.} {\it
Suppose that $c > |\om_2^{(0)}|$. Then:
 \begin{itemize}
\item[(i)] If $\ga \in \Ga$ satisfies $\ga(Y_c) \cap Y_c \not= \emptyset$ then $\ga \in \Ga \cap P_{r-1}$;
\item[(ii)] If $\ga \in \Ga(N)$ is such that $\ga(Y_c) \cap Y_c \not= \emptyset$ then $\ga$ has the shape
 $$\ga = \begin{array}{|c|c|}\hline
 1 & u_2,\ldots,u_r \\ \hline
 0 & \\
 \vdots & \ga' \\
 0 & \\ \hline
 \end{array}$$ 
where $u_2,\ldots,u_r \in NA$ and $\ga'$ runs through a finite subgroup of $\Ga'(N) = \Ga'\cap \Ga(N)$ consisting
of strictly upper triangular matrices (i.e., with ones on the diagonal). On the other hand, each $\ga$ of this form with
$\ga' = 1$ stabilizes $Y_c$. 
  \end{itemize}
}

 \begin{proof}
 \begin{itemize}
  \item[(i)] Let $\bo = (\om_1:\ldots:\om_r) \in Y_c$ be such that $\ga\bo = (\om'_1:\ldots:\om'_r) \in Y_c$ with 
  $\ga = (\ga_{i,j}) \in \Ga$ (recall that $\om_r = \om'_r = 1$). Let further $\La$ be the lattice $\La_{\bo} = \langle \om_1,\ldots,\om_r \rangle_A$ 
 and $\alpha:= {\rm aut}(\ga,\bo)$. Now $\{\om'_r,\ldots,\om'_1\}$ is a basis of $\alpha^{-1}\La$ and, since $\ga\bo \in Y_c$,
 $\om'_r,\ldots,\om'_2$ are the first $r-1$ elements of an SMB of $\alpha^{-1}\La$, so $\alpha \om'_r,\ldots,\alpha\om'_2$ are
 the first $r-1$ elements of an SMB of $\La$. Then $|\underset{1\leq j \leq r}{\sum} \ga_{i,j}\om_j| = |\alpha\om'_i| = |\om_i|$ holds
 for $i = 2,\ldots,r$ in view of (1.8.2). If $\ga_{i,1} \not= 0$ then $|\underset{1\leq j \leq r}{\sum}\ga_{i,j} \om_j| \geq |\om_1| > |\om_i|$, 
 contradiction.
  \item[(ii)] The entry $\ga_{1,1} = 1$ is obvious, as is the fact that each $\ga$ with $\ga' = 1$ stabilizes $Y_c$. The possible
  $\ga'$ are those that fix $X$. The stabilizer $\Ga'_X$ of $X$ in $\Ga'$ equals
   $$\{\ga' = (\ga_{i,j})_{2\leq i,j \leq r} \in \Ga' ~|~|\ga_{i,j}| \leq |\om_i^{(0)}/\om_j^{(0)}|\},$$
matrices with a block structure

$$\begin{BMAT}(b){|c:cc|}{|c:cc|}
     B_1 &  & \ast \\
      & B_2 &  \\
     0 &  & \ddots\\
  \end{BMAT}
\medskip $$

\noindent
 with zeroes below the blocks, and each block $B_k$ an invertible matrix over $\F$. The number of such $\ga'$ is finite, and the
 congruence condition $\ga' \equiv 1 (\bmod N)$ forces each block to equal 1. Hence $\ga'$ is strictly upper triangular.
  \end{itemize}
  \end{proof}
   \bigskip
   
(6.3) We write $G$ for the group that occurs in (6.2)(ii), i.e., $G:= \{\ga \in \Ga(N)~|~\ga(Y_c) \cap Y_c \not= \emptyset\}$, 
$G_1 := \{\ga \in G~|~\ga' = 1\}$ and $G' := G/G_1 =\Ga'_X$ for the group of possible $\ga'$.
 \bigskip
 
(6.4)  Consider the function $\bo \longmapsto t(\bo) := e^{-1}_{N\La}(\om_1)$ on $\Om^r$, where now $\La = \langle \om_2,\ldots,\om_r\rangle_A$. Its most important properties are:
 \bigskip
 
 (6.4.1) $t$ is well-defined, as $\om_1$ does not belong to $N\La$;
  \bigskip
  
 (6.4.2) it is holomorphic and has a unique strongly continuous extension to $\Om^r \cup \Om_{V_{r-1}}$, where $t \equiv 0$ on
 $\Om_{V_{r-1}}$;
  \bigskip
  
 (6.4.3) $t(\ga\bo) = t(\bo)$ for $\ga \in G$;
  \bigskip
  
 (6.4.4) Fix $\bo' = (\om_2:\ldots,\om_r) \in X$. Then the image of the map
  $$\begin{array}{rll}
   t_{\bo'}:\: \{\om \in C_{\infty}~|~(\om:\bo') \in Y_c\} & \lra & C_{\infty}\\
    \om & \longmapsto & t(\om:\bo')
     	\end{array}$$  
is a pointed ball $B_{\rho}^* := \{z \in C_{\infty}~|~0 < |z| \leq \rho\}$ for some $\rho = \rho(c) \in q^{\BQ}$,
and is independent of the choice of $\bo' \in X$. The function $c \longmapsto \rho(c)$ is strictly monotonically decreasing
with $\underset{c\to\infty}{\lim} \rho(c) = 0$.
 \bigskip
 
As for proofs, (6.4.1) is obvious, and (6.4.2) comes from trivial estimates. {\it Proof} of (6.4.3): As ${\rm aut}(\ga,\bo)=1$
for $\ga \in G$, $(\ga\bo)_2,\ldots,(\ga\bo)_r$ generate the same lattice $\La = \langle \om_2,\ldots,\om_r\rangle_A$. Hence
 $$t(\ga\bo) = e_{N\La}((\ga \bo)_1)^{-1} = e_{N\La}(\om_1)^{-1} = t(\bo),$$
 in view of $(\ga\bo)_1 \equiv \om_1(\bmod N\La)$ and the $N\La$-invariance of $e_{N\La}$. 
  \bigskip
  
 {\it Proof} of (6.4.4). The assertion that ${\rm im}(t_{\bo'})$ is a pointed ball $B_{\rho}^*$ for some $\rho$ is a general
fact of rigid analysis (see e.g. \cite{22} Lemma 10.9.1). Expansion of the product for $|e_{N\La}(\om)|$ shows that it depends
only on the distance $d(\om,K_{\infty}\La)$ and the values $|\om_2|,\ldots,|\om_r|$, as $\{\om_r,\ldots,\om_2\}$ is an SMB of
$\La$. Since the $|\om_i|$ are constant on $X$, the independence of $\rho(c)$ of the choice of $\bo'$ follows. The last statement is obvious. \hspace{2cm} $\Box$

  \bigskip
  
 {\bf 6.5 Proposition.} {\it 
 Let $\pi$ be the projection $\bo = (\om_1:\cdots:\om_r) \longmapsto \bo' = (\om_2:\ldots:\om_r)$ from $Y_c$ to $X$ and $\rho$ as in 
 $(6.4.4)$. Then $t \times \pi$ induces an isomorphism
  $$G\setminus Y_c \stackrel{\cong}{\lra} B_{\rho}^* \times (G'\setminus X)$$
 of analytic spaces.}
  \bigskip
  
\begin{proof}
For each $\bo' \in X$, $t_{\bo'}$ provides an isomorphism
 $$t_{\bo'}:\: (N\La)\setminus \{\om \in C_{\infty}~|~(\om:\bo') \in Y_c\} \stackrel{\cong}{\lra} B_{\rho}^*,$$
as it is bijective and 
 $$\frac{d}{d\om} t_{\bo'}(\om) = -t_{\bo'}(\om)^2 \not= 0 \quad \mbox{(since $\frac{d}{d\om} e_{N\La}(\om)=1$)}.
  \leqno{(6.5.1)}$$
  (Here $\La$ always denotes the lattice $\langle \om_2,\ldots,\om_r\rangle_A$ associated with $\bo$!) Therefore,
  also
   $$(t,\pi):\: G_1\setminus Y_c \stackrel{\cong}{\lra} B_{\rho}^* \times X$$
 is bijective, and is in fact an isomorphism, as, due to (6.5.1), its Jacobian matrix is invertible in each point.
 The group $G' = G/G_1$ acts on both sides (trivially on $B_{\rho}^*$) and due to (6.4.3), the map $(t,\pi)$ is 
 $G'$-equivariant. Therefore
  $$G\setminus Y_c = G'\setminus (G_1\setminus Y_c) \stackrel{\cong}{\lra} G' \setminus (B_{\rho}^* \times X) 
   = B_{\rho}^* \times (G'\setminus X).$$
 \end{proof}  
 
(6.6) As an admissible open affinoid in the smooth space $\Ga'(N)\setminus \Om^{r-1}$, the space $W:= G'\setminus X$ is
itself smooth and affinoid. Let $B = B_{\rho}$ be the unpunctured ball of radius $\rho$ and $Z := B \times W$.
 \bigskip
 
Eisenstein series extend uniquely (as strongly continuous functions and therefore, as $Z$ is smooth, by Bartenwerfer's 
criterion \cite{1} also as analytic functions) from $G\setminus Y_c \stackrel{\cong}{\lra} B^* \times W$ to $Z$. The
following result is a generalization of \cite{14} Korollar 2.2.
 \bigskip
 
{\bf 6.7 Proposition.} {\it 
Let $\fu = (u_1,\ldots,u_r) \in \MT(N)$, $u_1 = a/N$ with $\deg a < d = \deg N$. Then the Eisenstein series $E_{\fu} =
E_{1,\fu}$, regarded as a function on $Z$, has a zero of order $|a|^{r-1}$ along the divisor ($t=0$) of $Z$.
}   
 \bigskip
 
  \begin{proof}
  The proof has been given for the case $r=2$ in \cite{14}. There, equivalently, the pole order of $E_{\fu}(\bo)^{-1} =
  d_{\fu}(\bo)' = e_{\bo}(\fu \bo)$ has been determined, where the product expansion of $e_{\bo}$ was used. The proof
  generalizes without difficulty to the case of higher rank $r$. For another approach, see \cite{28}, Lemma 1.23.
   \end{proof}
   \bigskip
   
(6.8) Again by Bartenwerfer's criterion, the open embedding
 $$B^* \times W \stackrel{\cong}{\lra} G\setminus Y_c \hookrightarrow \Ga(N)\setminus \overline{\Om}^r$$
 extends to an injective map   
  $$i:\: Z = B \times W \hookrightarrow \Ga(N)\setminus \overline{\Om}^r,$$
  since $Z$ is smooth, thus normal. Its image is
   $${\rm im}(i) = G \setminus Y_c \cup G'\setminus X = G \setminus (Y_c \cup X).$$
We will show that $i$ is in fact an open embedding, i.e., an isomorphism of $Z=B \times W$ with ${\rm im}(i)$. That is,
given the class $[\bo'] \in W = G' \setminus X$ of $\bo' \in X$ with corresponding point $[\bo] = (0,[\bo']) \in Z$, we must
show that the canonical map
 $$\MO_{G\setminus(Y_c \cup X),[\bo]} \lra \MO_{Z,[\bo]}$$
of analytic local rings is an isomorphism. In fact, it suffices to show the corresponding isomorphism without dividing out
the group $G$. In other words, we must show that the canonical injection
 $$\MO_{Y_c \cup X,\bo} \hookrightarrow \MO_{B\times X,\bo} = \MO_{B,0} \widehat{\otimes} \MO_{X,\bo'}$$
is bijective. ($\widehat{\otimes}$ is  the topological tensor product of the two local rings, and we use ``$\bo$'' both for the
point $(\om_1 = 0, \bo')$ of $Y_c \cup X$ and the point $(t=0,\bo')$ of $B \times X$.) Now the left hand side contains
$\MO_{X,\bo'}$ via the projection $\pi:\: Y_c \cup X \lra X$, and it suffices to show that it also contains a uniformizer
$u$ of $Y_c\cup X$ along $X$, i.e., some $u$ that as a germ of a function on $B \times X$ presents a zero of order
$1$ along ($t=0$). Then $\MO_{Y_c\cup X,\bo}$ encompasses the local ring point ${\rm Sp}(C_{\infty}\langle u\rangle)$ at
$u=0$, that is, $\MO_{B,0}$, and we are done. By Proposition 6.7, we can take as $u$ the Eisenstein series $E_{\fu}$ with
$\fu = (\frac{1}{N},0,\ldots,0)\in \MT(N)$.
 \bigskip
 
Therefore we have shown the following result.
 \bigskip
 
{\bf 6.9 Theorem.} {\it  
Let $M_U(N) = \Ga(N)\setminus \Om_U$ with $U \in \mathfrak U_{r-1}$ be a cuspidal divisor on $\overline{M}^r(N)$, 
$\bo^{(0)}$ a point on $\Om_U$ with class $[\bo^{(0)}]$ in $M_U(N)$, and $\pi_U:\: \Om^r \lra \Om_U$ the canonical
projection. There exists an admissible open affinoid neighborhood $X$ of $\bo^{(0)}$ in $\Om_U$ and an admissible open
subspace $Y$ of $\Om^r$ characterized by
 $$Y = \{\bo \in \Om^r~|~\pi_U(\bo) \in X \mbox{ and $\bo$ sufficiently close to $\pi_U(\bo)$}\}$$
such that $Y \cup X$ is isomorphic with a product $B \times X$ and the subspace $Z := \Ga(N)\setminus (Y\cup X)$
of $\Ga(N)\setminus \overline{\Om}^r$ is isomorphic with $B\times W$, where $B$ is a ball and $W$ is an admissible
open affinoid neighborhood of $[\bo^{(0)}]$ in $M_U(N)$, quotient of $X$ by a finite group $G'$ of automorphisms.
The second projection $Z\lra W$ comes from $\pi_U$, divided out by the action of $\Ga(N)$, and the first projection
$Z\lra B$ is given by an explicit uniformizer.}
 \medskip
 
(If $U = V_{r-1}$ then $X,Y$ and $t$ are specified in (6.1) and (6.4), and 'sufficiently close' means 
$d(\om_1,\langle \om_2,\ldots,\om_r\rangle_{K_{\infty}}) \geq c$ for some constant $c$, or equivalently,
$|t(\bo)| \leq \rho$ for some $\rho$ depending on $c$.) \hspace{2cm}$\Box$
 \medskip
 
{\bf 6.10 Corollary.} {\it $[\bo^{(0)}]$ is a smooth point of $\overline{M}^r(N)$.}
 \medskip
 
\begin{proof}
$B$ and $W$ are smooth.
\end{proof} 
 \medskip
 
{\bf 6.11 Remark.} Theorem 6.9 allows to expand holomorphic functions on $\Om^r$ of weight $k$ for $\Ga(N)$ 
(so-called weak modular forms, see (7.1)) as Laurent series in $t$, where the coefficients are
holomorphic functions on $\Om_U$. This is crucial for the theory of modular forms. For the case $r=2$, see e.g.
\cite{24} or \cite{15}; for higher rank $r$, Basson and Breuer have started investigations in this direction in
\cite{2} and \cite{3}. See also \cite{4}, \cite{5}, \cite{6}.
 \medskip
 
{\bf 6.12 Remark.} Analogous tubular neighborhoods along cuspidal divisors may also be constructed for
$\overline{\widetilde{M}}^r(N)$. The proof for $\overline{M}^r(N)$ given in (6.4)--(6.8) may easily be adapted.   
 \vspace{0.5cm}
 
{\bf 7. Modular forms.}
 \vspace{0.3cm}
 
In this section, we define the ring of modular forms for $\Ga(N)$ and relate it with the ring of sections of the very ample
line bundle of $\overline{M}^r(N)$ given by the embedding $j_N$.
 \medskip
 
Theorem 7.9 gives several different descriptions of modular forms. The assumptions of the preceding sections remain
in force. Thus $r \geq 2$ and $N \in A$ is monic of degree $d \geq 1$. 
 \medskip
 
(7.1) Let $\MO(1)$ be the usual twisting line bundle on $\mathbb P = \mathbb P({\rm Eis}_1(N)^{\wedge})$ and
$\mathfrak M := j_N^*(\MO(1))$ its restriction to the subvariety $j_N:\: \overline{M}^r(N) \hookrightarrow \mathbb P$.
Tracing back the definitions, one sees that the sections of $\mathfrak M^{\otimes k}$ restricted to the open analytic
subspace $M^r(N) = \Ga(N)\setminus \Om^r$ are just the functions $f$ on $\Om^r$ subject to 
 \begin{itemize}
\item[(i)] $f$ is holomorphic;
\item[(ii)] for $\bo \in \Om^r$ and $\ga \in \Ga(N)$, the rule 
    $$f(\ga\bo) = {\rm aut}(\ga,\bo)^k f(\bo)$$
holds (${\rm aut}(\ga,\bo) = \underset{1\leq i \leq r}{\sum} \ga_{r,i}\om_i$, $\bo$ normalized such that $\om_r=1$).
 \end{itemize}  
For further use, we baptize functions $f$ satisfying (i) and (ii) as {\em weak modular forms} of weight $k$ for $\Ga(N)$.
Later we will define modular forms as weak modular forms that additionally satisfy certain boundary conditions discussed
below. 
 \medskip
 
(7.2) First, we recall the fundamental domain  $\BF$ for $\Ga$ on $\Om^r$:
$$\BF = \{\bo \in \Om^r~|~\{\om_r=1,\om_{r-1},\ldots,\om_1\}
\mbox{ is an SMB of its lattice $\La_{\bo}$}\}.\leqno{(7.2.1)}$$
The group $\Ga/\Ga(N)$ acts on $\overline{M}^r(N)$ and thus on weak modular forms of weight $k$ and level $N$
through
 $$f\longmapsto f_{[\ga]_k}, \mbox{ where } f_{[\ga]_k}(\bo) = {\rm aut}(\ga,\bo)^{-k} f(\ga \bo).\leqno{(7.2.2)}$$
(I.e., the formula is valid for $\ga \in \Ga$, and $\ga \in \Ga(N)$ acts trivially.) For simplicity, we choose and fix a 
system RS of representatives for $\Ga/\Ga(N)$.
 \medskip
 
(7.3) Let now $f$ be a weak modular form of some weight $k\in \N$. Consider the conditions on $f$:
 \begin{itemize}
\item[(a$^{\rm st}$)] $f$ extends to a holomorphic section of $\mathfrak M^{\otimes k}$ over $\overline{M}^r(N)$;
\item[(a)] $f$ extends to a strongly continuous section of $\mathfrak M^{\otimes k}$ over $\overline{M}^r(N)$ (that is, $f$
regarded as a $\Ga(N)$-invariant homogeneous function of weight $k$ on $\Psi^r$ has a strongly continuous extension
to $\overline{\Psi}^r$);
 \item[(b)] $f$ is integral over the ring ${\bf Mod} = C_{\infty}[g_1,\ldots,g_r]$ of modular forms of type 0 for $\Ga$;
  \item[(c)] $f$ along with all its conjugates $f_{[\ga]_k}$ ($\ga \in RS$) is bounded on $\BF$.
   \end{itemize}
 \medskip
 
{\bf 7.4 Theorem.} {\it   
For weak modular forms $f$ of weight $k$ and level $N$, we have the following implications:    
{\rm (a}$^{\rm st})$ $\Rightarrow$ {\rm (a)}   $\Leftrightarrow$ {\rm (b)} $\Leftrightarrow$ {\rm (c)}. If the variety $\overline{M}^r(N)$ happens to
be normal, then the four conditions are equivalent.}
 
  \begin{proof}
  \fbox{(a$^{\rm st}$) $\Rightarrow$ (a)} is trivial. 
  \medskip
  
  \fbox{(a) $\Rightarrow$ (b)}  The elementary symmetric functions in the $f_{[\ga]_k}$ ($\ga \in RS$) are invariant under
  $\Ga/\Ga(N)$, that is, under $\Ga$ and have a strongly continuous extension to $\overline{\Psi}^r$. As 
  $\Ga \setminus \overline{\Psi}^r = \overline{M}^r$ is normal, these extensions are in fact holomorphic by \cite{1} and
  therefore modular forms for $\Ga$. Hence $f$ satisfies an integral equation
   $$f^n+a_1f^{n-1}+\cdots + a_n = 0 \mbox{ with } a_i \in {\bf Mod}.\leqno{(7.4.1)}$$
   \fbox{(b) $\Rightarrow$ (c)} Suppose that $f$ is subject to an equation (7.4.1). Then it also holds for $f$ replaced
   by $f_{[\ga]_k}$, and $f$ along with all its conjugates is bounded on $\BF$, as the $a_i$ are (\cite{20} Proposition 1.8).
    \medskip
    
 \fbox{(c) $\Rightarrow$ (b)} As in \fbox{(a) $\Rightarrow$  (b)}, $f$ satisfies an equation of type (7.4.1) with the elementary
 symmetric functions in the $f_{[\ga]_k}$ ($\ga \in RS$) as coefficients $a_i$ up to sign. From the boundedness of the
 $f_{[\ga]_k}$ we conclude the boundedness of the $a_i$ on $\BF$, which in turn implies $a_i \in {\bf Mod}$ (\cite{20}
 Proposition 1.8).
  \medskip
  
 \fbox{(b) $\Rightarrow$ (a)} Suppose that $P(f) = 0$ with $P(X) = X^n+a_1X^{n-1} + \cdots + a_n$ and coefficients
 $a_1,\ldots,a_n \in {\bf Mod}$.  Regarding $f$ as a homogeneous and $\Ga(N)$-invariant function of weight $k$ on
 $\Psi^r$, we must show that it extends to a strongly continuous function on $\overline{\Psi}^r$.
  \medskip
  
Let $\bo = (U,i)$ be a boundary point, $\bo \in \overline{\Psi}^r \setminus \Psi^r$, and let $(\bo_{\ell})_{\ell \in \N}$ be
a sequence of elements of $\Psi^r$ that tends to $\bo$. Without restriction, replacing $f$ with some transform  
$f_{[\ga]_k}$ if necessary, we may assume that $\bo \in \Psi_{V_s}$ for some $1 \leq s < r$, $\bo = (0,\ldots,0,\om_{r-s+1},\ldots,
\om_r)$, where $\bo$ lies in the corresponding fundamental domain $\BF_s$, see (1.15.4). We will show:
 \begin{itemize}
  \item[(A)] $(f(\bo_{\ell}))_{\ell \in \N}$ converges in $C_{\infty}$.\\
  Therefore, $\overline{f}(\bo) := \underset{\ell\to\infty}{\lim} f(\bo_{\ell})= \underset{\underset{\bo'\to \bo}{\bo'\in \Psi^r}}{\lim}
  f(\bo')$ exists;
   \item[(B)] The so-defined extension $\overline{f}$ of $f$ to $\overline{\Psi}^r$ is strongly continuous, of weight $k$ and
   $\Ga(N)$-invariant.
    \end{itemize}
Put $\overline{a}_i := \underset{\ell\to \infty}{\lim} a_i(\bo_{\ell})$, which exists as $a_i$ is a modular form for $\Ga$, and let
$\overline{P}(X) = X^n+ \underset{1\leq i \leq n}{\sum} \overline{a}_iX^{n-i}$ be the limit polynomial. Elementary estimates
show that $f(\bo_{\ell})$ is close to a zero of $\overline{P}$ if $\ell \gg 0$. Hence the set of limit points of
$(f(\bo_{\ell}))_{\ell \in \N}$ is contained in the set $Z(\overline{P})$ of zeroes of $\overline{P}$, and each $f(\bo_{\ell})$ is close
to one of them for $\ell \gg 0$. Let $Y$ be a small neighborhood (w.r.t. the strong topology) of $\bo$ in 
$\overline{\Psi}^r$. Then for $Y$ and $\epsilon >0$ small enough,
   $$\Psi^r \cap Y = \underset{x\in Z(\overline{P})}{\stackrel{\bullet}{\bigcup}} Y_x,\leqno{(7.4.2)}$$
where $Y_x := \{\bo' \in \Psi^r \cap Y~|~|f(\bo')-x| < \epsilon\}$. We may further choose $Y$ such that    
 \medskip
$$\begin{array}{l}
 \Psi^r \cap Y = \{\bo' \in \Psi~|~(\om'_{r-s+1},\ldots,\om'_r) \\
 \mbox{lies in a  fixed connected open affinoid neighborhood $X$ of} \\
  \bo \mbox{ in } \Psi_{V_s} \mbox{ and }
   d(\om'_i,\langle\om'_{r-s+1},\ldots,\om'_r\rangle_{K_{\infty}}) \geq c \mbox{ for } 1 \leq i \leq r-s\}
   \end{array}\leqno{(7.4.3)}$$
 \noindent
for sufficiently large $c\in q^{\BQ}$. Such a set is connected as an analytic space.
 \bigskip 
 
 Now the occurrence of at least two different zeroes $x$ in (7.4.2) would contradict the connectedness of
 $\Psi^r \cap Y$. Hence there exists only one limit point $x$, which equals
  $$\overline{f}(\bo) := \lim_{\ell \to \infty} f(\bo_{\ell}) = \underset{\underset{\bo'\to \bo}{\bo'\in \Psi^r}}{\lim} f(\bo'),$$
and (A) is established. The fact (B) that $\overline{f}$  is strongly continuous is seen by a modification of the above argument, working
now with approximating sequences $(\bo_{\ell})_{\ell \in \N}$ for $\bo$ with $\bo_{\ell} \in \overline{\Psi}^r$. Also the
properties of weight $k$ and $\Ga(N)$-invariance turn over from $f$ to $\overline{f}$. 
 \medskip

Finally, suppose 
$\overline{M}^r(N)$ (and thus $\overline{\widetilde{M}}^r(N)$) is normal. Then the existence of a holomorphic extension
$\overline{f}$ of $f$ follows, again by Bartenwerfer's criterion, from the existence of a strongly continuous extension. Hence in
this case, (a) implies in fact (a$^{\rm st}$), and all four conditions are equivalent.
   \end{proof}
    \medskip
    
{\bf 7.5 Definition.} 
We define the Satake compactification $M^r(N)^{\rm Sat}$ of $M^r(N)$ as the normalization of 
$\overline{M}^r(N)$ in its function field $\MF_r(N)$. It is a normal projective $C_{\infty}$-variety provided with an embedding
$\iota:\: M^r(N) \hookrightarrow M^r(N)^{\rm Sat}$ and a finite birational morphism $\nu:\: M^r(N)^{\rm Sat} \lra \overline{M}^r(N)$
such that $\nu \circ \iota$ is the identity on $M^r(N)$. Likewise, we define $\widetilde{M}^r(N)^{\rm Sat}$ as the normalization
of $\overline{\widetilde{M}}^r(N)$ in its function field $\widetilde{\MF}_r(N)$. It has similar properties and is supplied with an action
of $\G_m$ such that $\G_m\setminus \widetilde{M}^r(N)^{\rm Sat} \stackrel{\cong}{\lra} M^r(N)^{\rm Sat}$.
 \medskip
 
{\bf 7.6 Corollary} (to the proof of Theorem 7.4; see also \cite{28}, proof of Proposition 1.23):
{\it The varieties $\overline{\widetilde{M}}^r(N)$ and $\overline{M}^r(N)$ are unibranched, that is, the canonical maps
 $$\widetilde{\nu}:\: \widetilde{M}^r(N)^{\rm Sat} \lra \overline{\widetilde{M}}^r(N) \mbox{ and } 
  \nu:\: M^r(N)^{\rm Sat} \lra \overline{M}^r(N)$$
are bijective.}
 \medskip
 
\begin{proof}
Since the open subspace $\widetilde{M}^r(N)$ of $\overline{\widetilde{M}}(N)$ is smooth, it suffices to consider boundary
points $[\bo]$ of $\overline{\widetilde{M}}(N)$ represented by $\bo$ as in the proof of \fbox{(b) $\Rightarrow$ (a)}. Then we have
to show that $[\bo]$ has at most one pre-image in the normalization $\widetilde{M}^r(N)^{\rm Sat}$. However, this follows
from the connectedness of the sets $\Psi^r \cap Y$ in (7.4.3). (If there were several pre-images of $[\bo]$ then 
$Y\setminus \{\bo\}$ and also $\Psi^r\cap Y$ had to split into several components for $Y$ sufficiently small.) The argument
for $\overline{M}^r(N)$ follows the same lines.
 \end{proof}
 \medskip
 
{\bf 7.7 Remark.} We know from (6.10) that the singular locus of $\overline{M}^r(N)$ is contained in
 $$\overline{M}^r_{\leq r-2}(N) := \bigcup_{U \in \mathfrak U \atop \dim U \leq r-2} M_U(N)$$ 
 and therefore has codimension $\geq 2$, as expected for a normal variety. Together with unibranchedness this suggests
 that $\overline{M}^r(N)$ should itself be normal, i.e., $M^r(N)^{\rm Sat} = \overline{M}^r(N)$. However, this is not a formal
 implication, and the question of normality of $\overline{M}^r(N)$ is still open.
  \medskip
  
(7.8) At least, $\nu:\: M^r(N)^{\rm Sat} \lra \overline{M}^r(N)$ is bijective by (7.6) and therefore (since it is a finite morphism)
a strong homeomorphism of the sets of $C_{\infty}$-points. As we don't know whether $\nu$ is always an isomorphism
(see Section 8 for examples), we make the following double definition.
 \medskip
 
(7.8.1) A {\em strong modular form} of weight $k$ and level $N$ is a weak modular form $f$ that satisfies condition
(7.3)(a$^{\rm st}$), that is, $f$ extends to a holomorphic section of $\mathfrak M^{\otimes k}$ over $\overline{M}^r(N)$.
A {\em modular form} of weight $k$ and level $N$ is a weak modular form $f$ that satisfies (7.3)\,(a), i.e., the boundary
condition is relaxed to: $f$ extends to a strongly continuous section of $\mathfrak M^{\otimes k}$, or equivalently 
(by (7.6) and Bartenwerfer's criterion), $f$ extends holomorphically to a section of $\nu^*(\mathfrak M^{\otimes k})$ over 
the Satake compactification $M^r(N)^{\rm Sat}$. 
 \medskip
 
(7.8.2) We let ${\bf Mod}_k´^{\rm st}(N)$ resp. ${\bf Mod}_k(N)$ be the  $C_{\infty}$-spaces of (strong) modular forms
of weight $k$ and 
 $${\bf Mod}^{\rm st}(N) = \underset{k\geq 0}{\oplus} {\bf Mod}^{\rm st}_k(N),\: {\bf Mod}(N) = \underset{k\geq 0}{\oplus} {\bf Mod}_k(N)$$
the corresponding graded $C_{\infty}$-algebras. Then 
 $${\bf Eis}(N) \subset {\bf Mod}^{\rm st}(N) \subset {\bf Mod}(N).\leqno{(7.8.3)}$$
 The common quotient field of the three rings is the field $\widetilde{\MF}_r(N)$. By (7.4) the following criterion holds.
  \medskip

{\bf 7.9 Theorem.} {\it Let $f$ be a weak modular form of weight $k$ and level $N$. Then the following three conditions are
equivalent:
 \begin{itemize}
  \item[(a)] $f \in {\bf Mod}(N)$;
  \item[(b)] $f$ is integral over ${\bf Mod} = C_{\infty}[g_1,\ldots,g_r]$;
  \item[(c)] $f$ and all its conjugates $f_{[\ga]_k}$ ($\ga \in RS$) are bounded on the fundamental domain $\BF$.
  \end{itemize}
  }
  \medskip
  
(7.10) Let $J$ be the ideal of Eisenstein relations in $R = {\rm Sym}({\rm Eis}_1(N))$, see (5.4). To the exact sequence
 $$0 \lra J \lra R \lra {\bf Eis}(N) \lra 0$$
corresponds an exact sequence of sheaves (in the algebraic sense) on the variety $\mathbb P = {\rm Proj}(R)$
 $$0 \lra \mathfrak J \lra \MO_{\BP} \lra \MO_{\overline{M}^r(N)} \lra 0,\leqno{(7.10.1)}$$
where we regard the structure sheaf $\MO_{\overline{M}^r(N)} $ of $\overline{M}^r(N)$ as a sheaf on $\BP$ with support
in $\overline{M}^r(N) \hookrightarrow \BP$. It remains exact upon tensoring with the sheaf $\MO(k)= \MO(1)^{\otimes k}$ over $\BP$,
where $k>0$. As $\MO(1)$ restricted to $\overline{M}^r(N)$ is the sheaf $\mathfrak M$ of strong modular forms, we find the
exact sequence
 $$0 \lra \mathfrak J(k) \lra  \MO_{\BP}(k) \lra \mathfrak M(k) \lra 0.\leqno{(7.10.2)}$$
The first part of its exact cohomology sequence reads: 
 $$\begin{array}{c}
 0 \lra H^0(\BP,\mathfrak J(k)) \lra H^0(\BP,\MO_{\BP}(k)) \stackrel{\alpha}{\lra} H^0(\BP,\mathfrak M(k)) \lra\vspace{0.2cm}\\
 \lra H^1(\BP,\mathfrak J(k)) \lra H^1(\BP,\MO_{\BP}(k)) \lra \ldots
  \end{array} $$
  Now, 
    \begin{itemize}
   \item $H^1(\BP,\MO_ {\BP}(k))$ vanishes (see e.g. \cite{25} III Theorem 5.1);
 \item $H^0(\BP,\mathfrak M(k)) = H^0(\overline{M}^r(N),\mathfrak M^{\otimes k}) = {\bf Mod}^{\rm st}_k(N)$;
 \item ${\rm im}(\alpha)$ is the subspace ${\bf Eis}_k(N)$ of strong modular forms that belong to the Eisenstein algebra
 ${\bf Eis}(N)$. 
    \end{itemize}
 Hence $H^1(\BP,\mathfrak J(k))$ measures the difference between ${\bf Eis}_k(N)$ and ${\bf Mod}^{\rm st}(N)$. By standard
 properties (\cite{25} III Theorem 5.2), $H^1(\BP,\mathfrak J(k))$ vanishes for large $k$. As it is always finite-dimensional, we see:
  \medskip
  
 {\bf 7.11 Corollary.} {\it 
 For $k$ sufficiently large, ${\bf Mod}_k^{\rm st}(N)$ agrees with its subspace ${\bf Eis}_k(N)$. In particular, the 
 Eisenstein algebra ${\bf Eis}(N)$ has finite codimension in the algebra ${\bf Mod}^{\rm st}(N)$. }
  \medskip
  
 {\bf 7.12 Corollary.} {\it
 $\overline{M}^r(N)={\rm Proj}({\bf Eis}(N)) = {\rm Proj}({\bf Mod}^{\rm st}(N))$.}
  \medskip
  
 \begin{proof}
 The first equality has been shown in Section 5, the second is a formal consequence of the definition of ${\bf Mod}^{\rm st}(N)$, 
 but follows also from $\dim({\bf Mod}^{\rm st}(N)/{\bf Eis}(N)) < \infty$.
  \end{proof}   
   \medskip
   
(7.13) Consider the projective variety ${\rm Proj}({\bf Mod}(N))$ attached to the algebra of modular forms. It is normal
(as ${\bf Mod}(N)$ is integrally closed), provided with a natural map to ${\rm Proj}({\bf Mod}^{\rm st}(N)) = \overline{M}^r(N)$, and birational
with $\overline{M}^r(N)$, and thus agrees with the normalization, i.e., 
 $$M^r(N)^{\rm Sat} = {\rm Proj}({\bf Mod}(N)). \leqno{(7.13.1)}$$   
 
{\bf 7.14 Corollary.}{\it 
The three assertions are equivalent:
 \begin{itemize}
\item[(i)] $\overline{M}^r(N)$ is normal.
\item[(ii)] ${\bf Mod}^{\rm st}(N) = {\bf Mod}(N)$;
\item[(iii)] ${\bf Mod}^{\rm st}(N)$ has finite codimension in ${\bf Mod}(N)$.
  \end{itemize}
  }
 \medskip
 
 \begin{proof}
 (i) $\Rightarrow$ (ii) is the last assertion of Theorem 7.4, and (ii) $\Rightarrow $ (iii) is trivial. Suppose that (iii) holds.
 Then ${\rm Proj}({\bf Mod}(N)) = {\rm Proj}({\bf Mod}^{\rm st}(N)) = \overline{M}^r(N)$, and it follows that the latter is
 normal.	
  \end{proof}  
   \medskip
   
{\bf 7.15 Corollary.} {\it  
Suppose that $\overline{M}^r(N)$ fails to be normal. Then there exist arbitrarily large weights $k$ such that ${\bf Mod}_k(N)$ is
strictly larger than its subspace ${\rm Eis}_k(N)$.}
 \medskip
 
 \begin{proof}
As all the ${\bf Mod}_k(N)$ have finite dimension, this follows from the last corollary.
 \end{proof}   
  \medskip
  
(7.16) We conclude this section with an observation about ${\bf Mod}^{\rm st}(N)$. Given $\fu \in \MT(N)$, we let
${\bf Eis}(N)_{E_{\fu}}$ be the localization w.r.t. $E_{\fu}$, i.e., ${\rm Spec}({\bf Eis}(N)_{E_{\fu}})$ is the open subvariety
of $\overline{\widetilde{M}}^r(N)$ where $E_{\fu}$ doesn't vanish. Hence
 $$\overline{\widetilde{M}}^r(N) = \bigcup_{\fu \in \MT(N)} {\rm Spec}({\bf Eis}(N)_{E_{\fu}}) =\bigcup_{\fu \in N^{-1}S}
  \ldots,$$
where $S$ is the set of representatives of $(A/N)^r_{\rm prim}/\BF^*$ given in (3.5). Similarly, 
 $$\overline{M}^r(N) = \bigcup_{\fu \in N^{-1}S} \overline{M}^r(N)_{(E_{\fu} \not= 0)}.$$
 A weak modular form of weight $k$ extends to a section of $\mathfrak M^{\otimes k}$ if and only if its restriction to each
 $\overline{M}^r(N)_{(E_{\fu} \not= 0)}$ has the corresponding property, that is, belongs to ${\bf Eis}(N)_{E_{\fu}}$. Therefore
 we may describe ${\bf Mod}^{\rm st}(N)$ as the intersection
  $${\bf Mod}^{\rm st}(N) = \bigcap_{\fu \in N^{-1}S} {\bf Eis}(N)_{E_{\fu}} \leqno{(7.16.1)}$$
in $\widetilde{\MF}(N)$.
 \vspace{0.5cm}
 
{\bf 8. Examples and concluding remarks.}
 \vspace{0.3cm}
 
The preceding immediately raises a number of important questions and desiderata. 
 \medskip
 
{\bf 8.1 Question.} Do the Eisenstein and Satake compactifications $\overline{M}^r(N)$ and $M^r(N)^{\rm Sat}$ always
coincide, i.e., is $\overline{M}^r(N)$ always normal?
 \medskip
 
(8.2) Describe the singularities of both compactifications and construct natural desingularizations together with a modular
interpretation! (See \cite{31} for some results.)
 \medskip
 
(8.3) How far do the algebras ${\bf Eis}(N)$, ${\bf Mod}^{\rm st}(N)$, ${\bf Mod}(N)$ differ, if at all? Describe their 
Hilbert functions, that is, the dimensions of their pieces in dimension $k$, and find presentations for these algebras!
 \medskip
 
Almost nothing about these questions is known when the rank $r$ is larger than 2. We will briefly present the 
state-of-the-art in the case where \fbox{$r=2$}, which we assume until (8.9). Here the $\overline{M}^2(N)$ are smooth curves
\cite{11}, so the Eisenstein and Satake compactification agree, and therefore ${\bf Mod}^{\rm st}(N) = {\bf Mod}(N)$. The genera of
the $\overline{M}^2(N)$ have been determined by Goss \cite{23} and, with a different method, by the author \cite{13}.
 \medskip
 
(8.4) Let $d \geq 1$ be the degree of $N$, and suppose that
 $$N = \prod_{1 \leq i \leq t} \mathfrak p_i^{s_i}$$
is the prime decomposition, where the $\mathfrak p_{i}$ are different monic prime polynomials. As in (3.6), write
$q_i = q^{\deg\,\mathfrak p_i}$. We define
 $$\la(N) := \prod_{1\leq i \leq t}q_i^{2s_i-2} (q_i^2-1),\leqno{(8.4.1)}$$
which appears in the formulas below. (Note that $\la(N) = \varphi(N)\epsilon(N)$ with the arithmetic functions $\varphi$, $\epsilon$ 
defined in \cite{17} 1.5.) Then the numbers $g(N) = $ genus of the modular curve 
$\overline{M}^2(N)$, $c(N)=c_2(N)=$ number
of cusps of $M^2(N)$, $\deg(\mathfrak M) =$ degree of the line bundle of modular forms over $\overline{M}^2(N)$ and 
$\dim\,{\bf Mod}_k(N) = \dim\,{\bf Mod}_k^{\rm st}(N)$ are given by
 \medskip
 
(8.4.2) $g(N) = 1+\la(N)(q^d-q-1)/(q^2-1)$;
 \medskip
 
(8.4.3) $c(N) = \la(N)/(q-1)$;
 \medskip
 
(8.4.4) $\deg(\mathfrak M) = \la(N)q^d/(q^2-1)$;
 \medskip
 
(8.4.5) $\dim\,{\bf Mod}_k(N) = ((k-1)q^d+q+1)\la(N)/(q^2-1)$.
 \medskip
 
Here (8.4.2) and (8.4.3) may be found in \cite{23} and \cite{13} (such data for other Drinfeld modular curves are collected 
in \cite{17}) and (8.4.4) is from \cite{15} VII 6.1. The last formula (8.4.5) is an immediate consequence of the Riemann-Roch 
theorem provided that $k \geq 2$; for $k=1$, Riemann-Roch and Serre duality yield only
 $$\dim\,{\bf Mod}_1(N) =c(N)+\dim\,H^1(\overline{M}^2(N),\mathfrak M),\leqno{(8.4.6)}$$    
where $c(N) = \dim\,{\bf Eis}_1(N) = \dim\,{\rm Eis}_1(N)$ and $\dim\,H^1(\overline{M}^2(N),\mathfrak M) =
\dim\,{\bf Mod}_1^2(N)$ with the space ${\bf Mod}_1^2(N)$ of double cuspidal (double zeroes at the cusps) modular forms
of weight 1 (\cite{15} p. 92). However by the next result, (8.4.5) is valid for $k=1$, too.
 \medskip
 
{\bf 8.5 Proposition.} {\it  
For $r=2$ we have ${\rm Eis}_1(N) = {\bf Mod}_1(N)$, of dimension $c(N)$.}
 \medskip
 
As a proof of this basic fact so far has not been published, we give a brief sketch here.
 \medskip
 
\begin{proof}
By Corollary 4.7, $\dim\,{\rm Eis}_1(N) = c(N)$. Since moreover the space of cusp forms ${\bf Mod}_1^1(N)$ is a complement
of ${\rm Eis}_1(N)$ in ${\bf Mod}_1(N)$ (this is a consequence of Proposition 4.6), it suffices to show that there are
no non-trivial cusp forms of weight 1.
 \medskip
 
Assume that $f \in {\bf Mod}_1^1(N)$. Then $f^p$ is a cusp form of weight $p \geq 2$ for $\Ga(N)$, where $p = {\rm char}(\BF)$.
For $0 \leq i \leq p-2$, the residue ${\rm res}_e \om^{i} f^p(\om) d\om$ of the differential form $\om^{i}f^p(\om)d\om$ on 
$\Om = \Om^2$ at the oriented edge $e$ of the Bruhat-Tits tree of ${\rm PGL}(2,K_{\infty})$ vanishes for each $e$, as is
immediate from the definition of ${\rm res}_e$, see \cite{34} Definition 9. Hence the image ${\rm res}(f^p)$ under
Teitelbaum's isomorphism ({\it loc.\,cit.}, Theorem 16) of cusp forms of weight $p$ with the space of cocyles of a certain
type vanishes, and so do $f^p$ and $f$ itself. 
 \end{proof}   
  \medskip
  
{\bf 8.6 Remark.}  All the formulas and results in (8.4) and (8.5) have generalizations
 \begin{itemize}
  \item to other congruence subgroups of $\Ga$, e.g., Hecke congruence subgroups $\Ga_0(N)$, $\Ga_1(N)$, etc., see \cite{17};
   \item to more general Drinfeld rings $A$ than $A=\F[T]$, e.g., $A$ the affine ring of an elliptic curve over $\F$, see \cite{15}
   pp. 92--93.
  \end{itemize}
As to the relationship between ${\bf Eis}(N)$ and ${\bf Mod}(N)$, there is the following result of Cornelissen \cite{9}.
 \medskip
 
{\bf 8.7 Theorem} (Cornelissen). {\it Let $r=2$. The algebra ${\bf Mod}(N)$ of modular forms for $\Ga(N)$ is generated by 
${\bf Mod}_1(N) = {\rm Eis}_1(N)$ and the space ${\bf Mod}_2^1(N)$ of cusp forms of weight 2.}
 \medskip
 
In fact, it is not difficult (using \cite{14} Korollar 2.2) to show that ${\bf Mod}_2^1(N)$ above may be replaced with the
space ${\bf Mod}_2^2(N)$ of  double cuspidal forms of weight 2, which under $f(\om) \longmapsto f(\om)d\om$
corresponds to the $g(N)$-dimensional space of holomorphic differentials on $\overline{M}^2(N)$. Still, this doesn't
completely answer the question (8.3) of whether ${\bf Eis}(N) = {\bf Mod}(N)$ in this case. The only positive results in this
direction seem to be the following two examples.
 \medskip
 
{\bf 8.8 Example.} Suppose that $r=2$ and $d = \deg\,N=1$. Then $g(N)=0$, that is, $\overline{M}^2(N)$ is a projective
line, and $\deg(\mathfrak M)=q$. Therefore, ${\bf Eis}(N)={\bf Mod}(N)$ and $\dim\,{\bf Mod}_k(N) = 1+kq$ in this case.
From this, a presentation of ${\bf Mod}(N)$ may be derived (\cite{9}, see also \cite{35}, \cite{36}, which also study the 
spaces ${\bf Mod}_k(N)$ as modules under the action of $\Ga/\Ga(N) = {\rm GL}(2,\F)$).
 \medskip
 
{\bf 8.9 Example} (Cornelissen \cite{8}). Let again $r=2$, $q=2$, and $d = \deg\,N=2$. Then ${\bf Eis}(N) = {\bf Mod}(N)$.
 \medskip
 
As the possible genera $g(N)$ here are positive (they may take the values 4, 5 and 6), this case is less trivial than (8.8).
The equality of the two rings, i.e., the projective normality of the Eisenstein embedding $\overline{M}^2 \hookrightarrow \BP$,
is based on a numerical criterion of Castelnuovo. Unfortunately, the validity of this argument is strictly limited to the
requirements of Example 8.9.   
 \medskip
 
For the next example, we return to the general case, where $r \geq 2$ is arbitrary.
 \medskip
 
{\bf 8.10 Example.} Suppose that $d = \deg\,N=1$. After a coordinate change, we may assume $N=T$. This case has been
extensively studied by Pink and Schieder \cite{32}. Actually, they consider a certain $\F$-variety $Q_V$, but which after base 
extension with $C_{\infty}$ and some translational work may also be seen as our $\overline{M}^r(T)$. Their results
(overlapping with (8.8) if $r=2$) give
 \begin{itemize}
  \item ${\bf Eis}(T) = {\bf Mod}(T)$ (Theorem 1.7 in \cite{32}), i.e., the normality of $R_V = {\bf Eis}(T)$;
  \item  a presentation through generators and relations (Theorem 1.6);
  \item the Hilbert function of ${\bf Eis}(T)$ (Theorem 1.10).
  \end{itemize}
Further, they construct and discuss a desingularization $B_V$ of $Q_V$ (Section 10). Hence, concerning our questions
(8.1)--(8.3), nothing is left to desire. But note that these satisfactory and complete results refer only to the (isolated?)
case where $\deg\,N = 1$. None of our current knowledge excludes the possibility that always ${\bf Eis}(N) = {\bf Mod}(N)$, or 
the sheer opposite possibility that the two rings differ almost always and (8.8), (8.9), (8.10) are just extreme boundary cases.

  Ernst-Ulrich Gekeler\\ Fachrichtung Mathematik der Universit\"at des Saarlandes\\ Campus E2 4, 66123 Saarbr\"ucken, Germany,
  gekeler@math.uni-sb.de

 \end{document}